\documentclass[number,seceqn,dvips]{arxbj}
\usepackage{cursive,mathrsfs}  %uplgreek}

\aid{0}
\volume{16}
\issue{4}
\pubyear{2010}
\firstpage{1262}
\lastpage{1293}
\doi{10.3150/10-BEJ258}

\makeatletter
\newtheorem{theo}{Theorem}
\newtheorem{lem}[theo]{Lemma}
\newtheorem{prop}[theo]{Proposition}
\newremark{rem}[theo]{Remark}

\newcommand{\N}{\mathbb{N}}
\newcommand{\R}{\mathbb{R}}
\newcommand{\C}{\mathbb{C}}
\newcommand{\bP}{\mathbb{P}}
\def\ra{\rangle}
\def\la{\langle}
\def\HH{\mathfrak{H}}
\def\e{\varepsilon}
\def\1{{\mathbf{1}}}
\def\sfrac#1#2{#1/#2}
\def\afrac#1#2{(#1)/#2}

\makeatother

\begin{document}
\begin{frontmatter}

\title{Limit theorems for nonlinear functionals of Volterra processes via white noise analysis}
\runtitle{Limit theorems for nonlinear functionals of Volterra processes via white noise analysis}

\begin{aug}
\author[a]{\fnms{S\'ebastien} \snm{Darses}\corref{}\thanksref{a}\ead[label=e1]{darses@cmi.univ-mrs.fr}},
\author[b]{\fnms{Ivan} \snm{Nourdin}\thanksref{b}\ead[label=e2]{ivan.nourdin@upmc.fr}}\and
\author[c]{\fnms{David} \snm{Nualart}\thanksref{c}\ead[label=e3]{nualart@math.ku.edu}}
\runauthor{S. Darses, I. Nourdin and D. Nualart}
\address[a]{Universit\'e Aix-Marseille I, 39 rue Joliot Curie, 13453 Marseille Cedex 13, France.\\\printead{e1}}
\address[b]{Laboratoire de Probabilit\'es et Mod\`eles Al\'eatoires, Universit\'e Pierre et Marie Curie (Paris VI),
Bo\^ite Courrier 188, 4 place Jussieu, 75252 Paris Cedex 05, France. \printead{e2}}
\address[c]{Department of Mathematics, University of Kansas, Lawrence, KS 66045, USA.\\\printead{e3}}
\pdfauthor{Sebastien Darses, Ivan Nourdin, David Nualart}
\end{aug}

% HISTORY:
\received{\smonth{4} \syear{2009}}

% ABSTRACT
\begin{abstract}
By means of white noise analysis, we prove some limit theorems for nonlinear functionals
of a given Volterra process. In particular, our results apply to fractional Brownian motion (fBm) and should be
compared with the classical convergence results of the 1980s due to Breuer, Dobrushin, Giraitis, Major,
Surgailis and Taqqu, as well as the recent advances concerning the construction of a
L\'{e}vy area for fBm
due to Coutin, Qian and Unterberger.
\end{abstract}

% KEYWORDS

\end{frontmatter}

%s1 ###
\section{Introduction}\label{intro}

Fix $T>0$ and let $B=(B_t)_{t\geq 0}$ be a fractional Brownian motion
with Hurst index $H\in(0,1)$,
defined on some probability space $(\Omega,\mathcal{B},P)$. Assume that
$\mathcal{B}$ is the completed $\sigma$-field generated by $B$.
Fix an integer $k\geq 2$ and, for $\e>0$, consider
%e1.1 ###
\begin{equation}\label{geps}
G_\e = \e^{-k(1-H)}\int_0^T h_k
\biggl(\frac{B_{u+\e}-B_u}{\e^{H}}\biggr)\,\mathrm{d}u.
\end{equation}
Here, and in the rest of this paper,
%e1.2 ###
\begin{equation}\label{herm-pol}
h_k(x)=(-1)^k \mathrm{e}^{x^2/2}\frac{\mathrm{d}^k}{\mathrm{d}x^k} (\mathrm{e}^{-x^2/2} )
\end{equation}
stands for the $k$th Hermite polynomial. We have
$h_2(x)=x^2-1$, $h_3(x)=x^3-3x$ and so on.

Since the seminal works \cite{BM,DM,GS,taqqu75,taqqu79} by Breuer, Dobrushin, Giraitis, Major, Surgailis and
Taqqu, the  following three convergence results are classical:
\begin{itemize}
\item if $H<1- \frac 1{2k}$, then
%e1.3 ###
\begin{equation}\label{cv<}
\bigl((B_t)_{t\in[0,T]}
,\e^{k(1-H) -1/2} G_{\e}  \bigr)
\displaystyle\mathop{\stackrel{\mathrm{Law}}{\longrightarrow}}_{\e\to 0}
\bigl((B_t)_{t\in[0,T]}, N  \bigr),
\end{equation}
where $N\sim\mathscr{N} (0,T
\times  k!\int_0^T\rho^k(x)\,\mathrm{d}x )$ is
independent of $B$, with $\rho(x)=\frac12 (|x+1|^{2H}+|x-1|^{2H}-2|x|^{2H} )$;
\item if $H=1- \frac 1{2k}$, then
%e1.4 ###
\begin{equation}\label{cv=}
\biggl((B_t)_{t\in[0,T]}
, \frac{G_{\e}}
{
\sqrt{\log(1/\e)}
}  \biggr) \displaystyle\mathop{\stackrel{  \mathrm{Law}}{\longrightarrow}}_{\e\to 0  }
\bigl((B_t)_{t\in[0,T]}, N  \bigr),
\end{equation}
where $N\sim\mathscr{N} (0,T
\times 2k!(1-\frac1{2k})^k(1-\frac1k)^k )$ is independent of $B$;
\item if  $H>1-\frac1{2k}$, then
%e1.5 ###
\begin{equation}\label{cv>}
G_{\e} \displaystyle\mathop{\stackrel{L^2(\Omega)}{\longrightarrow}}_{\e\to 0}
Z_T^{(k)},
\end{equation}
where $Z^{(k)}_T$ denotes the Hermite random variable of order $k$; see Section~\ref{sec31}
for its definition.
\end{itemize}

Combining (\ref{cv<}) with the fact that $\sup_{0<\e\leq 1}E [|\e^{k(1-H)-1/2}G_\e|^p ]<\infty$ for all $p\geq 1$
(use  the boundedness of $\operatorname{Var}(\e^{k(1-H)-1/2}G_\e)$ and a
classical hypercontractivity argument), we have, for all $\eta\in
L^2(\Omega )$ and if $H<1-\frac1{2k}$,
that
\[
\e^{k(1-H)-1/2}E[\eta G_\e]
\mathop{\longrightarrow}_{\e\to 0}
E(\eta  N)=E(\eta)E(N)=0
\]
(a similar statement holds in the critical case $H=1-\frac1{2k}$).
This means that $\e^{k(1-H)-1/2} G_\e$ converges \textit{weakly}
in $L^2(\Omega)$ to zero. The following question then arises.  Is
there a normalization of $G_\e$ ensuring that it converges \textit{
weakly} towards a \textit{non-zero} limit when $H\leq 1-\frac1{2k}$? If
so, then what can be said
about the limit? The first purpose  of the present paper is to
provide an answer to this question in the framework of \textit{white
noise analysis}.

In \cite{nualartwhite}, it is shown that for all $H\in(0,1)$, the
time derivative $\dot{B}$ (called the \textit{fractional white noise})
is a distribution in the sense of Hida.  We also refer to Bender
\cite{bender}, Biagini \textit{et al.}~\cite{BOSW} and references therein for
further works on the fractional white noise.

Since we have $E(B_{u+\e}-B_u)^2=\e^{2H}$, observe that $G_\e$
defined in (\ref{geps}) can be rewritten as
%e1.6 ###
\begin{equation}
G_{\e}=\int_0^T  \biggl(\frac{B_{u+\e}-B_u}{\e} \biggr)^{\diamond k}\,\mathrm{d}u,  \label{e1bis}
\end{equation}
where $(\ldots)^{\diamond k}$ denotes the $k$th Wick product. In
Proposition~\ref{thm-fbm} below, we will show that for all
$H\in (\frac12-\frac1k,1 )$,
%e1.7 ###
\begin{equation}\label{cvwick}
\lim_{\e\to 0}\int_0^T
\biggl(\frac{B_{u+\e}-B_u}{\e} \biggr)^{\diamond k}\,\mathrm{d}u
=\int_0^T \dot{B}_u^{\diamond k}\,\mathrm{d}u,
\end{equation}
where the limit is in the $(\mathcal{S})^*$ sense.

In particular, we observe two different types of asymptotic results for
$G_{\e}$ when $H\in (\frac12-\frac1k,1- \frac 1{2k} )$: convergence (\ref{cvwick}) in $(\mathcal{S})^*$  to
a Hida distribution,
and convergence (\ref{cv<}) in law to a normal law, with rate  $\e^{ 1/2 -  k(1-H)  }$.
On the other hand, when $H\in (1-\frac1{2k},1 )$, we obtain from (\ref{cv>}) that the
Hida distribution
$\int_0^T \dot{B}_s^{\diamond k}\,\mathrm{d}s$ turns out to be the square-integrable random variable
$Z_T ^{(k)}$, which is  an interesting result in its own right.

In Proposition~\ref{CVS*}, the convergence (\ref{cvwick}) in $(\mathcal{S})^*$ is proved for a general class of Volterra processes of the
form
%e1.8 ###
\begin{equation}\label{def-vol}
\int_0^t K(t,s)\, \mathrm{d}W_s,\qquad  t\geq 0,
\end{equation}
where   $W$ stands for a standard Brownian motion, provided the kernel $K$
satisfies some suitable conditions; see Section~\ref{volterra}.

We also  provide a new proof of the convergence (\ref{cv<})
based on the recent general criterion for the convergence in
distribution to a normal law of a sequence of multiple stochastic
integrals established by Nualart and Peccati \cite{NP} and by
Peccati and Tudor \cite{PT}, which avoids the classical method of
moments.

In two recent papers \cite{MR1,MR2},    Marcus and Rosen have
obtained central and non-central limit theorems for a functional of
the form (\ref{geps}), where $B$ is a mean zero Gaussian process
with stationary increments
such that the covariance function of $B$, defined by
$\sigma^2(|t-s|)=\operatorname{Var}(B_t-B_s)$, is either convex (plus some
additional regularity conditions), concave or given by
$\sigma^2(h)=h^r$ with $1<r<2$. These theorems include the
convergence  (\ref{cv<}) and, unlike our simple proof, are  based on
the method of moments.

In the second part of the paper, we develop a similar analysis for
functionals of two independent fractional Brownian motions (or, more
generally, Volterra processes) related to the L\'evy area. More
precisely, consider two \textit{independent} fractional Brownian
motions $B^{(1)}$ and $B^{(2)}$ with (for simplicity) the same Hurst
index $H\in(0,1)$. We are interested in the convergence, as $\e\to
0$, of
%e1.9 ###
\begin{equation}\label{tilde}
\widetilde{G}_\e:=\int_0^T B^{(1)}_u \frac{B^{(2)}_{u+\e}-B^{(2)}_u}{\e}\, \mathrm{d}u
\end{equation}
and
%e1.10 ###
\begin{equation}\label{breve}
\breve{G}_\e:=\int_0^T  \biggl(\int_0^u \frac{B^{(1)}_{v+\e}-B^{(1)}_v}{\e}\, \mathrm{d}v \biggr)
\frac{B^{(2)}_{u+\e}-B^{(2)}_u}{\e} \,\mathrm{d}u .
\end{equation}
Note that $\widetilde{G}_\e$ coincides with the $\e$-integral
associated with the forward Russo--Vallois integral $\int_0^T B^{(1)}\,\mathrm{d}^-B^{(2)}$;
see, for example, \cite{RVLN} and references therein.
Over the last decade, the convergences of $\widetilde{G}_\e$
and $\breve{G}_\e$ (or of related families of random variables) have
been intensively studied.
Since $\e^{-1}\int_0^u (B^{(1)}_{v+\e}-B^{(1)}_v)\,\mathrm{d}v$ converges pointwise to $B^{(1)}_u$
for any $u$, we could think that the asymptotic behaviors of $\widetilde{G}_\e$ and $\breve{G}_\e$
are very close as $\e\to 0$. Surprisingly, this is not the case. Actually, only the result for
$\breve{G}_\e$ agrees with the seminal result of Coutin and Qian \cite{CQ} (that is,
we have convergence of $\breve{G}_\e$ in $L^2(\Omega)$ if and only if $H>1/4$)
and with the recent result by Unterberger \cite{unterbergertcl}
(that is, adequately renormalized,
$\breve{G}_\e$ converges in law if $H<1/4$). More precisely:
\begin{itemize}
\item if $H<1/4$, then
%e1.11 ###
\begin{equation}\label{star10}
\bigl( \bigl(B^{(1)}_t,B^{(2)}_t\bigr) _{t\in [0,T]},    \e^{1/2-2H} \breve{G}_{\e}  \bigr)
\displaystyle\mathop{\stackrel{  \mathrm{Law}}{\longrightarrow}}_{\e\to 0}\hspace*{1pt}
\bigl( \bigl(B^{(1)}_t,B^{(2)}_t\bigr) _{t\in [0,T]},  N \bigr),
\end{equation}
where    $N  \sim  \mathscr{N} (0,T\breve{\sigma}^2_H )$ is independent of
$(B^{(1)},B^{(2)})$ and
\begin{eqnarray*}
&&\breve{\sigma}^2_H=\frac{1}{4(2H+1)(2H+2)}
\int_\R (|x+1|^{2H}+|x-1|^{2H}-2|x|^{2H} )\\
&&\hspace*{119pt}{}\times (2|x|^{2H+2}-|x+1|^{2H+2}-|x-1|^{2H+2} )\,\mathrm{d}x;
\end{eqnarray*}
\item if $H=1/4$, then
%e1.12 ###
\begin{equation}\label{star10crit}
\biggl( \bigl(B^{(1)}_t,B^{(2)}_t\bigr) _{t\in [0,T]},
\frac{ \breve{G}_{\e}}{\sqrt{\log(1/\e)}}    \biggr)
\displaystyle\mathop{\stackrel{  \mathrm{Law}}{\longrightarrow}}_{\e\to 0}\hspace*{1.5pt}
\bigl( \bigl(B^{(1)}_t,B^{(2)}_t\bigr) _{t\in [0,T]},  N \bigr),
\end{equation}
where    $N  \sim  \mathscr{N} (0,T/8 )$ is independent of
$(B^{(1)},B^{(2)})$;
\item if $H> 1/4$, then
%e1.13 ###
\begin{equation}\label{star11}
\breve{G}_{\e} \displaystyle\mathop{\stackrel{  L^2(\Omega)  }{\longrightarrow}}_{\e\to 0}
\int_0^T  B^{(1)}_u\diamond \dot{B}^{(2)}_u \,\mathrm{d}u=\int_0^T B^{(1)}_u\,\mathrm{d}B^{(2)}_u;
\end{equation}
\item for all $H\in(0,1)$, we have
%e1.14 ###
\begin{equation}\label{star12}
\breve{G}_{\e} \displaystyle\mathop{\stackrel{  (\mathcal{S})^*}{\longrightarrow}}_{\e\to 0}
\int_0^T  B^{(1)}_u\diamond \dot{B}^{(2)}_u \,\mathrm{d}u.
\end{equation}
\end{itemize}
However, for $\widetilde{G}_\e$, we have, in contrast:
\begin{itemize}
\item if $H<1/2$, then
%e1.15 ###
\begin{equation}\label{star1}
\bigl( \bigl(B^{(1)}_t,B^{(2)}_t\bigr) _{t\in [0,T]},  \e^{1/2-H} \widetilde{G}_{\e}   \bigr)
\displaystyle\mathop{\stackrel{  \mathrm{Law}}{\longrightarrow}}_{\e\to 0  }\hspace*{1.5pt}
\bigl( \bigl(B^{(1)}_t,B^{(2)}_t\bigr) _{t\in [0,T]},  N\times S \bigr),
\end{equation}
where
\[
S=\sqrt{\int_0^\infty (|x+1|^{2H}+|x-1|^{2H}-2|x|^{2H} )\,\mathrm{d}x\times\int_0^T \bigl(B^{(1)}_u\bigr)^2\,\mathrm{d}u}
\]
and
$N\sim\mathscr{N}(0,1)$, independent of $(B^{(1)},B^{(2)})$;
\item if $H\geq 1/2$, then
%e1.16 ###
\begin{equation}\label{star3}
\widetilde{G}_{\e}
\displaystyle\mathop{\stackrel{  L^2(\Omega)  }{\longrightarrow}}_{\e\to 0}
\int_0^T  B^{(1)}_u\diamond \dot{B}^{(2)}_u \,\mathrm{d}u=\int_0^T B^{(1)}_u\,\mathrm{d}B^{(2)}_u;
\end{equation}
\item for all $H\in(0,1)$, we have
%e1.17 ###
\begin{equation}\label{star4}
\widetilde{G}_{\e} \displaystyle\mathop{\stackrel{  (\mathcal{S})^*}{\longrightarrow}}_{\e\to 0}
\int_0^T  B^{(1)}_u\diamond \dot{B}^{(2)}_u \,\mathrm{d}u.
\end{equation}
\end{itemize}

Finally, we study the convergence, as $\e\to 0$, of the so-called $\e$\textit{-covariation} (following
the terminology of \cite{RVLN})
defined by
%e1.18 ###
\begin{equation}
\widehat{G}_\e:=
\int_0^T \frac{B^{(1)}_{u+\e}-B^{(1)}_u}{\e}\times\frac{B^{(2)}_{u+\e}-B^{(2)}_u}{\e} \,\mathrm{d}u\label{hat}
\end{equation}
and we get:
\begin{itemize}
\item if $H<3/4$, then
%e1.19 ###
\begin{equation}\label{star5}
\bigl( \bigl(B^{(1)}_t,B^{(2)}_t\bigr) _{t\in [0,T]}, \e^{3/2-2H} \widehat{G}_{\e}   \bigr)
\displaystyle\mathop{\stackrel{  \mathrm{Law}}{\longrightarrow}}_{\e\to 0  }\hspace*{1pt}
\bigl( \bigl(B^{(1)}_t,B^{(2)}_t\bigr) _{t\in [0,T]}, N  \bigr)
\end{equation}
with  $N \sim \mathscr{N}(0,T\widehat{\sigma}^2_H)$ independent of $(B^{(1)},B^{(2)})$  and
\[
\widehat{\sigma}^2_H=\frac14
\int_\R (|x+1|^{2H}+|x-1|^{2H}-2|x|^{2H} )^2\,\mathrm{d}x;
\]
\item if $H=3/4$, then
%e1.20 ###
\begin{equation}\label{star6}
\biggl( \bigl(B^{(1)}_t,B^{(2)}_t\bigr) _{t\in [0,T]},\frac{\widehat{G}_{\e}}
{ \sqrt{\log(1/\e)} }    \biggr)  \displaystyle\mathop{\stackrel{  \mathrm{Law}}{\longrightarrow}}_{\e\to 0  }\hspace*{1.5pt}
\bigl( \bigl(B^{(1)}_t,B^{(2)}_t\bigr) _{t\in [0,T]}, N  \bigr)\\
\end{equation}
with  $N \sim \mathscr{N}(0,9T/32)$ independent of $(B^{(1)},B^{(2)})$;
\item if $H>3/4$, then
%e1.21 ###
\begin{equation}\label{star7}
\widehat{G}_{\e} \displaystyle\mathop{\stackrel{  L^2(\Omega)  }{\longrightarrow}}_{\e\to 0}
\int_0^T  \dot{B}^{(1)}_u\diamond \dot{B}^{(2)}_u \,\mathrm{d}u;
\end{equation}
\item for all $H\in(0,1)$, we have
%e1.22 ###
\begin{equation}\label{star8}
\widehat{G}_{\e} \displaystyle\mathop{\stackrel{  (\mathcal{S})^* }{\longrightarrow}}_{\e\to 0}
\int_0^T  \dot{B}^{(1)}_u\diamond \dot{B}^{(2)}_u \,\mathrm{d}u.
\end{equation}
\end{itemize}

The paper is organized as follows. In Section~\ref{sec2},
we introduce some preliminaries on white noise analysis.
Section~\ref{volterra} is devoted to the study,  using the language and  tools of the
previous section, of the asymptotic behaviors of $G_\e$, $\widetilde{G}_\e$
and $\widehat{G}_\e$
in the (more general) context where $B$ is a Volterra process.
Section~\ref{sec4} is concerned with the fractional Brownian motion case.
In Section~\ref{sec5} (resp., Section~\ref{sec6}), we prove (\ref{cv<}) and (\ref{cv=})
(resp., (\ref{star10}), (\ref{star10crit}), (\ref{star1}), (\ref{star5}) and (\ref{star6})).

%s2 ###
\section{White noise functionals}\label{sec2}

In this  section, we present some preliminaries on white noise analysis.
The classical approach to white noise distribution theory is
to endow the space of tempered
distributions $\mathcal{S}'(\R)$ with a~Gaussian measure $\bP$ such that, for  any rapidly decreasing function $\eta\in \mathcal{S}(\R)$,
\[
\int_{\mathcal{S}'(\R)}\mathrm{e}^{\mathrm{i}\la x,\eta\ra}\bP(\mathrm{d}x)=\mathrm{e}^{-\sfrac{|\eta|_0^2}{2}}.
\]
Here, $|\cdot|_0$ denotes the norm in  $L^2(\R)$
and  $\langle  \cdot,  \cdot \rangle$ the dual
pairing between $\mathcal{S}'(\R)$  and  $\mathcal{S}(\R)$.
The existence of such a measure is ensured by Minlos'
theorem
\cite{kuo}.

In this way, we can consider the probability space $(\Omega,\mathcal{B},\bP)$, where $\Omega=\mathcal{S}'(\R)$.
The pairing $\langle x, \xi \rangle$ can be extended, using the norm of  $L^2(\Omega)$, to any function
$\xi \in L^2(\R)$. Then, $W_t=\la\cdot,\1_{[0,t]}\ra$ is a two-sided Brownian motion
(with the convention that $ \1_{[0,t]}= -\1_{[t,0]}$ if $t<0$) and
for any   $\xi \in L^2(\R)$,
\[
\la\cdot,\xi \ra= \int_{-\infty}^\infty \xi\, \mathrm{d}W=I_1(\xi)
\]
is the Wiener integral of $\xi$.

Let $\Phi\in L^2(\Omega)$.  The classical Wiener chaos expansion of $\Phi$ says that  there exists a
sequence of  symmetric square-integrable  functions
$\phi_n \in L^2(\R^n)$ such that
%e2.1 ###
\begin{equation}   \label{wiener}
\Phi= \sum_{n=0} ^\infty  I_n(\phi_n),
\end{equation}
where  $I_n$ denotes the multiple stochastic integral.

%s2.1 ###
\subsection{The space of Hida distributions}
Let us recall some basic facts concerning tempered distributions.
Let $(\xi_n)_{n=0}^\infty$ be the orthonormal basis of $L^2(\R )$ formed by the Hermite functions
given by
%e2.2 ###
\begin{equation}\label{xi_n}
\xi_n(x)=\curpi^{-1/4}(2^nn!)^{-1/2}\mathrm{e}^{-x^2/2}h_n(x), \qquad x\in\R,
\end{equation}
where $h_n$ are the Hermite polynomials defined in (\ref{herm-pol}).
The following two facts can immediately be checked: (a)
there exists a constant $K_1>0$ such that $\|\xi_n\|_\infty\leq K_1(n+1)^{-1/12}$;
(b) since $\xi'_n=\sqrt{\frac{n}2}\xi_{n-1}-\sqrt{\frac{n+1}{2}}\xi_{n+1}$,
there exists a constant $K_2>0$ such that $\|\xi'_n\|_\infty\leq K_2 n^{5/12}$.

Consider the positive self-adjoint operator $A$ (whose inverse is Hilbert--Schmidt) given by
$A=-\frac{\mathrm{d}^2}{\mathrm{d}x^2}+(1+x^2)$. We have  $A\xi_n=(2n+2)\xi_n$.

For any $p\ge 0$, define the space $\mathcal{S}_p(\R)$ to be the domain of the closure of $A^p$.
Endowed with the norm $|\xi|_p:=|A^p\xi|_0$, it is a Hilbert space.
Note that the norm $|\cdot|_p$ can be expressed as follows, if one uses the orthonormal basis $(\xi_n)$:
\[
|\xi|_p^2=\sum_{n=0}^\infty  \langle \xi,\xi_n\rangle^2(2n+2)^{2p}.
\]
We denote by $\mathcal{S}'_p(\R)$ the dual of $\mathcal{S}_p(\R)$. The norm in  $S'_p(\R)$
is given by (see  \cite{ob}, Lemma 1.2.8)
\[
|\xi|^2_{-p}=\sum_{n=0} ^\infty |\la \xi,A^{-p}\xi_n\ra |^2
=\sum_{n=0}^\infty  \langle \xi,\xi_n\rangle^2(2n+2)^{-2p}
\]
for any $\xi\in \mathcal{S}'_p(\R)$.
One can show that the projective limit of the spaces $\mathcal{S}_{p}(\R)$, $p\geq 0$, is  $\mathcal{S}(\R)$, that the inductive
limit of the spaces $ \mathcal{S}_p(\R)'$, $p\ge 0$, is $\mathcal{S}'(\R)$ and that
\[
\mathcal{S}(\R)\subset L^2(\R)\subset \mathcal{S}'(\R)
\]
is a Gel'fand triple.

We can now introduce the Gel'fand triple
\[
(\mathcal{S})\subset L^2(\Omega)\subset(\mathcal{S})^*,
\]
via the second quantization operator $\Gamma(A)$. This is an unbounded  and densely defined operator
on $L^2(\Omega)$ given by
\[
\Gamma(A)\Phi= \sum_{n=0} ^\infty  I_n(A^{\otimes n}\phi_n),
\]
where $\Phi$ has the Wiener chaos expansion (\ref{wiener}).
If $p\geq 0$, then we denote by $(\mathcal{S})_{p}$ the space of random variables
$\Phi \in L^{2}(\Omega)$ with Wiener chaos expansion (\ref{wiener})  such that
\[
\| \Phi  \| ^p_{p}:=E [ |  \Gamma(A)^p \Phi |^2 ] =
\sum_{n=0} ^\infty n! |\phi_n|_p^2 <\infty.
\]
In the above formula,  $|\phi_n|_p$ denotes the norm in  $\mathcal{S}_p(\R) ^{\otimes n}$.
The projective limit of the spaces $(\mathcal{S})_{p}$, $p\geq 0$, is called the space
of test functions and is denoted by $(\mathcal{S})$. The inductive limit of the spaces $(\mathcal{S})_{-p}$, $p\geq 0$, is
called the space of Hida distributions and is denoted by $(\mathcal{S})^{\ast }$.
The elements of $(\mathcal{S})^*$ are called \textit{Hida distributions}.
The main example is the time derivative of the Brownian motion, defined as  $\dot{W}_t=\langle \cdot,\delta_t\rangle$.
One can show that  $|\delta_t|_{-p}<\infty$ for some $p>0$.

We denote by  $\langle\!\langle\Phi, \Psi\rangle\!\rangle$ the dual  pairing associated with the spaces $(\mathcal{S})$ and  $(\mathcal{S})^*$.
On the other hand (see \cite{ob}, Theorem 3.1.6),   for any $\Phi\in (\mathcal{S})^*$, there exist $\phi_n\in  \mathcal{S}(\R^n)'$ such that
\[
\langle\!\langle \Phi,\Psi \rangle\!\rangle=\sum_{n=0}^\infty n!\langle \phi_n,\psi_n\rangle,
\]
where $\Psi=\sum_{n=0}^\infty  I_n(\psi_n)\in (\mathcal{S})$. Moreover, there exists $p>0$ such that
\[
\|\Phi\|_{-p}^2=\sum_{n=0}^\infty  n!|\phi_n|_{-p}^2.
\]
Then, with a convenient abuse of notation,  we  say that $\Phi$ has a generalized Wiener chaos expansion
of the form  (\ref{wiener}).

%s2.2 ###
\subsection{The $S$-transform}
A useful tool to characterize elements in $(\mathcal{S})^*$ is the $S$-transform.
The Wick exponential of a Wiener integral $I_1(\eta)$, $\eta \in L^2(\R)$, is defined by
\[
: \mathrm{e}^{I_1(\eta)} :\ = \mathrm{e}^{ I_1(\eta) - |\eta|_0^2/2}.
\]
The  $S$-transform of an element $\Phi\in (\mathcal{S})^*$ is then defined by
\[
S(\Phi)(\xi)=\big\langle\big\langle\Phi,:\mathrm{e}^{I_1(\xi)}:\big\rangle\big\rangle,
\]
where   $\xi\in \mathcal{S}(\R)$.
One can easily see that the $S$-transform is injective on $(\mathcal{S})^*$.

If $\Phi\in L^2(\Omega)$, then $S(\Phi)(\xi)=E[\Phi:\mathrm{e}^{I_1(\xi)}:]$. For instance, the $S$-transform of the Wick exponential is
\[
S\bigl(:\mathrm{e}^{I_1(\eta)}:\bigr)(\xi)=\mathrm{e}^{\langle\eta,\xi\rangle}.
\]
Also,  $S(W_t)(\xi)=\int_0^t\xi(s)\,\mathrm{d}s$ and $S(\dot{W}_t)(\xi)=\xi(t)$.

Suppose that  $\Phi\in (\mathcal{S})^*$ has a generalized Wiener chaos expansion of the form
(\ref{wiener}).  Then,   for any  $\xi\in \mathcal{S}(\R)$,
\[
S(\Phi)(\xi)= \sum_{n=0}^\infty  \la \phi_n, \xi^{\otimes n} \ra,
\]
where the series converges absolutely (see   \cite{ob}, Lemma 3.3.5).

The Wick product of two functionals $\Psi=\sum_{n=0}^\infty I_n(\psi_n)$ and
$\Phi=\sum_{n=0}^\infty I_n(\phi_n)$ belonging to $(\mathcal{S})^*$ is defined  as
\[
\Psi\diamond\Phi =\sum_{n,m=0}^\infty I_{n+m}(\psi_n \otimes \phi_m).
\]
It can be proven that  $\Psi\diamond\Phi  \in (\mathcal{S})^*$. The following is an important property of the
$S$-transform:
%e2.3 ###
\begin{equation}\label{prod}
S(\Phi \diamond \Phi)(\xi) =S(\Phi)(\xi)S(\Psi)(\xi).
\end{equation}
If $\Psi$, $\Phi$ and $\Psi\diamond\Phi $ belong to  $L^2(\Omega)$, then we   have  $E[\Psi\diamond\Phi]=E[\Psi] E[\Phi]$.

The following is a useful characterization theorem.
%t1
\begin{theo}\label{F}
A function $F$ is the $S$-transform of an element $\Phi\in (\mathcal{S})^*$ if and only if the following conditions are satisfied:
\begin{enumerate}[(1)]
\item[(1)] for any $\xi,\eta\in \mathcal{S}$, $z\mapsto F(z\xi+\eta)$ is holomorphic on $\C$;
\item[(2)] there exist non-negative numbers $K,a$ and $p$ such that for all $\xi\in \mathcal{S}$,
\[
|F(\xi)|\leq K \exp(a|\xi|_p^{2}).
\]
\end{enumerate}
\end{theo}

\begin{pf}
See \cite{kuo}, Theorems 8.2 and  8.10.
\end{pf}

In order to study the convergence of a sequence in $(\mathcal{S})^*$, we can use its $S$-transform, by virtue of the following theorem.
%t2
\begin{theo}\label{th_kuo}
Let $\Phi_n \in (\mathcal{S})^*$ and $S_n=S(\Phi_n)$. Then, $\Phi_n$ converges in $(\mathcal{S})^*$ if and only if the following conditions are satisfied:
\begin{enumerate}
\item[(1)] $\lim_{n\to \infty}S_n(\xi)$ exists for each $\xi\in \mathcal{S}$;
\item[(2)] there exist non-negative numbers $K,a$ and $p$ such that for all $n\in \N$, $\xi\in (\mathcal{S})$,
\[
|S_n(\xi)|\leq K \exp(a|\xi|_p^{2}).
\]
\end{enumerate}
\end{theo}

\begin{pf}
See \cite{kuo}, Theorem 8.6.
\end{pf}

%s3 ###
\section{Limit theorems for Volterra processes}\label{volterra}

%s3.1 ###
\subsection{One-dimensional case}

Consider a Volterra process $B=( B_t )_{t\geq 0}$ of the form
%e3.1 ###
\begin{equation}\label{star}
B_t=\int_0^t K(t,s)\,\mathrm{d}W_s,
\end{equation}
where $K(t,s)$ satisfies $\int_0^tK(t,s)^2 \,\mathrm{d}s<\infty$ for all $t>0$ and $W$
is the Brownian motion defined on the white noise probability space introduced in the last section.
Note that the $S$-transform of the random variable $B_t$ is given by
%e3.2 ###
\begin{equation} \label{e1aa}
S(B_t)(\xi)=\int_0^{t}K(t,s)\xi(s)\,\mathrm{d}s
\end{equation}
for any $\xi\in  \mathcal{S}(\R)$.
We introduce the following assumptions on the kernel $K$:
\begin{enumerate}[$(\mathrm{H}_{\mathrm{1}})$]
\item[$(\mathrm{H}_{\mathrm{1}})$] $K$ is continuously differentiable on $\{0<s<t<\infty \}$ and for any $ t>0$,  we have
\[
\int_0^t     \bigg|\frac{\partial K}{\partial t}(t,s) \bigg|(t-s)\,\mathrm{d}s <\infty;
\]
\item[$(\mathrm{H}_{\mathrm{2}})$] $k(t)=\int_0^{t}K(t,s)\,\mathrm{d}s $ is  continuously differentiable on $(0,\infty)$.
\end{enumerate}
Consider the  operator $K_+$ defined by
\[
K_+\xi(t)=k'(t)\xi(t)+\int_0^t \frac{\partial K}{\partial t}(t,r)\bigl(\xi(r)-\xi(t)\bigr)\,\mathrm{d}r,
\]
where $t>0$ and  $\xi\in  \mathcal{S}(\R)$.  From Theorem~\ref{F},  it follows that  the linear
mapping $\xi\rightarrow K_+\xi(t)$ is the $S$-transform of a Hida distribution.
More precisely,
according to \cite{nualartwhite}, define the function
%e3.3 ###
\begin{equation}  \label{C}
C(t)=|k'(t)|+\int_0^t \bigg|\frac{\partial K}{\partial t}(t,r) \bigg|(t-r)\,\mathrm{d}r, \qquad t\geq 0,
\end{equation}
and observe that the following estimates hold (recall
the definition (\ref{xi_n}) of $\xi_n$):\vspace*{2pt}
%e3.4 ###
\begin{eqnarray}\label{e11}
|K_+\xi (t)|
& \leq & C(t) (\|\xi\|_\infty+\|\xi'\|_\infty)  \nonumber \\[1pt]
& \leq & C(t) \sum_{n=0}^\infty|\la\xi,\xi_n\ra| (\|\xi_n\|_\infty+\|\xi'_n\|_\infty) \nonumber\\[1pt]
&\leq& C(t)M\sum_{n=0}^\infty|\la \xi,\xi_n\ra|(n+1)^{5/12}\\[1pt]
& \leq & C(t) M   \sqrt{ \sum_{n=0}^\infty|\la\xi,\xi_n\ra|   ^2 (2n+2) ^{17/6}   }\sqrt{\sum_{n=0}^\infty(n+1)^{-2}}\nonumber\\[1pt]
&= & C(t) M   |\xi|_{17/12}\nonumber
\end{eqnarray}
for some constants  $M>0$ whose values are not always the same from one line to the next.

We have the following preliminary  result.
%l3
\begin{lem}\label{def-int}
Fix an integer $k\geq 1$. Let $B$ be a Volterra process with kernel $K$ satisfying the conditions
$(\mathrm{H_1})$ and $(\mathrm{H_2})$. Assume, moreover, that $C$ defined by (\ref{C})
belongs to $L^k([0,T])$.
The function $\xi\mapsto \int_0^T (K_+\xi(s))^k\,\mathrm{d}s$ is then the $S$-transform of an
element of $(\mathcal{S})^*$. This element is denoted by $\int_0^T \dot{B}_u^{\diamond k}\,\mathrm{d}u$.
\end{lem}
\begin{pf}
We use Theorem~\ref{F}. Condition (1) therein is immediately checked, while for condition~(2), we just write,
using (\ref{e11}),\vspace*{2pt}
\[
\bigg| \int_0^T (K_+\xi(s))^k\, \mathrm{d}s \bigg|\leq \int_0^T | K_+\xi(s)|^k\, \mathrm{d}s\leq M|\xi|_{17/12}\int_0^T C^k(s)\,\mathrm{d}s.
\]
\upqed
\end{pf}

Fix an integer $k\ge 1$ and consider  the following, additional, condition.
\begin{enumerate}[$(\mathrm{H}_{\mathrm{3}}^{k})$]
\item[$(\mathrm{H}_{\mathrm{3}}^{k})$]  The maximal function $D(t)=\sup_{0<\e\leq \e_0}\frac1\e\int_t^{t+\e}C(s)\,\mathrm{d}s$ belongs
to $L^k([0,T])$ for any $T>0$ and for some $\e_0>0$.
\end{enumerate}
We can now state the main result of this section.
%p4
\begin{prop}\label{CVS*}
Fix an integer $k\geq 1$. Let $B$ be a Volterra process with kernel $K$ satisfying the conditions
$\mathrm{(H_1)}$, $\mathrm{(H_2)}$ and $(\mathrm{H}_3^{k})$. The following convergence then holds:\vspace*{2pt}
\[
\int_0^T
\biggl(\frac{B_{u+\e}-B_u}{\e} \biggr)^{\diamond k}\,\mathrm{d}u \displaystyle\mathop{\stackrel{  (\mathcal{S})^*}{\longrightarrow}}_{\e\to 0}
\int_0^T \dot{B}_u^{\diamond k}\,\mathrm{d}u.
\]
\end{prop}

\begin{pf}
Fix $\xi\in{\mathcal{S}}(\R)$ and set
\[
S_\e(\xi)= S \biggl(\int_0^T  \biggl(\frac{B_{u+\e}-B_u}{\e} \biggr)^{\diamond k}\,\mathrm{d}u \biggr)(\xi).
\]
From linearity and property (\ref{prod}) of the $S$-transform, we obtain
%e3.5 ###
\begin{equation}  \label{e8}
S_\e(\xi)=\int_0^T \frac{(S(B_{u+\e}-B_u)(\xi))^{k}}{\e^k}\,\mathrm{d}u.
\end{equation}
Equation (\ref{e1aa}) yields
%e3.6 ###
\begin{equation}\label{e314}
S(B_{u+\e}-B_u)(\xi)= \int_0^{u+\e} K(u+\e,r) \xi(r)\,\mathrm{d}r -\int_0^{u} K(u,r) \xi(r)\,\mathrm{d}r.
\end{equation}
We claim that
%e3.7 ###
\begin{equation}  \label{e7}
\int_0^{u+\e} K(u+\e,r) \xi(r)\,\mathrm{d}r -\int_0^{u} K(u,r) \xi(r)\,\mathrm{d}r
=\int_u^{u+\e}K_+\xi(s) \,\mathrm{d}s.
\end{equation}
Indeed, we can write
\begin{eqnarray}\label{e6}
\hspace*{-10pt}\nonumber \int_u^{u+\e}K_+\xi(s) \,\mathrm{d}s
&=&\int_u^{u+\e} k'(s) \xi(s) \,\mathrm{d}s
+   \int_u^{u+\e}\biggl(  \int_0^s   \frac{\partial K}{\partial s}(s,r)\bigl(\xi(r)-\xi(s)\bigr) \,\mathrm{d}r  \biggr) \,\mathrm{d}s \\[-8pt]\\[-8pt]
&=& A^{(1)}_u +A^{(2)}_u.\nonumber
\end{eqnarray}
We have, using Fubini's theorem, that
\begin{eqnarray}
A^{(2)}_u&=& - \int_u^{u+\e}  \mathrm{d}s \nonumber
\int_0^s \mathrm{d}r \,\frac{\partial K}{\partial s}(s,r)     \int_r^s   \mathrm{d}\theta \,  \xi'(\theta)     \\[-8pt]\\[-8pt]
&=& -\int_0^{u+\e} \mathrm{d}\theta\, \xi'(\theta)\int_0^\theta
\mathrm{d}r\,
\bigl(  K(u+\e, r)- K(\theta \vee u,r)  \bigr).\nonumber
\end{eqnarray}
This can be rewritten as
\begin{eqnarray}\label{g1}
A^{(2)}_u&=& -\int_0^u   \bigl( K(u+\e,r) -K(u,r)  \bigr) \bigl(\xi(u)-\xi(r)\bigr) \,\mathrm{d}r  \nonumber \\[-8pt]\\[-8pt]
&&{}-\int_u^{u+\e}  \mathrm{d}\theta\, \xi'(\theta)  \int_0^\theta \mathrm{d}r  \, \bigl( K(u+\e,r) -K(\theta,r)  \bigr).\nonumber
\end{eqnarray}
On the other hand, integration by parts yields
\begin{eqnarray}  \label{e3}
A^{(1)}_u&=&\xi(u+\e) \int_0^{u+\e}  K(u+\e,r) \,\mathrm{d}r \nonumber\\[-8pt]\\[-8pt]
&&{} -\xi(u) \int_0^{u}  K(u,r) \,\mathrm{d}r
-\int_u^{u+\e}\,\mathrm{d}s\, \xi'(s) \int_0^s \,\mathrm{d}r \, K(s,r).\nonumber
\end{eqnarray}
Therefore, adding (\ref{e3}) and (\ref{g1}) yields
%e3.8 ###
\begin{eqnarray}\label{g3}
A^{(1)}_u +A^{(2)}_u
&=& \xi(u+\e) \int_0^{u+\e} K(u+\e,r) \,\mathrm{d}r
-  \xi(u ) \int_0^{u } K(u,r) \,\mathrm{d}r \nonumber \\
&&{}    -\int_0^u   \bigl( K(u+\e,r) -K(u,r)  \bigr) \bigl(\xi(u)-\xi(r)\bigr) \,\mathrm{d}r   \\
&&{} -\int_u^{u+\e}  \mathrm{d}\theta\, \xi'(\theta) \int_0^\theta     K(u+\e,r)  \,\mathrm{d}r.\nonumber
\end{eqnarray}
Note that, by integrating by parts, we have
%e3.9 ###
\begin{eqnarray}\label{g2}
&&-\int_u^{u+\e}  \mathrm{d}\theta\, \xi'(\theta) \int_0^\theta     K(u+\e,r)\,   \mathrm{d}r\nonumber
\\
&&\quad
=-\xi(u+\e) \int_0^{u+\e} K(u+\e,r)\, \mathrm{d}r
+ \xi(u ) \int_0^{u } K(u+\e,r) \, \mathrm{d}r
\\
&&\qquad{}+ \int_u^{u+\e} K(u+\e,r) \xi(r) \, \mathrm{d}r.\nonumber
\end{eqnarray}
Thus, substituting  (\ref{g2})  into (\ref{g3}), we obtain
\[
A^{(1)}_u +A^{(2)}_u  = \int_0^{u+\e} K(u+\e,r) \xi(r) \, \mathrm{d}r-  \int_0^{u } K(u,r) \xi(r) \, \mathrm{d}r,
\]
which completes the proof of  (\ref{e7}).
As a consequence, from  (\ref{e8})--(\ref{e7}),  we obtain
\[
S_\e (\xi)=\int_0^{T}  \biggl( \frac 1 \e \int_u^{u+\e} K_+\xi (s)\, \mathrm{d}s   \biggr)^k\, \mathrm{d}u.
\]
On the other hand, using (\ref{e11}) and the definition of the maximal function $D$, we get
%e3.10 ###
\begin{eqnarray}\label{v12}
\sup_{0< \e\leq \e_0}   \bigg| \frac 1 \e \int_u^{u+\e} K_+\xi(s)\, \mathrm{d}s \bigg |^k
&\leq& M^k |\xi|^k_{17/12} \sup_{0<\e\leq \e_0} \biggl(\frac1\e\int_u^{u+\e}C(s)\, \mathrm{d}s \biggr)^k\nonumber\\[-8pt]\\[-8pt]
&=&M^k|\xi|^k_{17/12}D^k(u).\nonumber
\end{eqnarray}
Therefore, using hypothesis   $(\mathrm{H}_3^{k})$  and the dominated convergence theorem, we have
%e3.11 ###
\begin{equation}
\lim_{\e\to0} S_\e(\xi)=\int_0^{T} (K_+\xi (s))^k\,\mathrm{d}s.  \label{e12}
\end{equation}
Moreover,\vspace*{-2pt} since $|S_\e(\xi)|\leq M^k|\xi|^k_{17/12}\int_0^T D^k(u)\,\mathrm{d}u$ for all $0<\e\leq \e_0$ (see
(\ref{v12})), conditions (1) and (2) in Proposition~\ref{CVS*} are fulfilled. Consequently,
$\e^{-k}\int_0^T (B_{u+\e}-B_u)^{\diamond k}\,\mathrm{d}u$ converges in $({\mathcal{S}}^*)$ as
$\e\to 0$.

To complete the proof, it suffices to observe that the right-hand side of (\ref{e12}) is, by definition (see Lemma~\ref{def-int}), the
$S$-transform of $\int_0^T \dot{B}_s^{\diamond k}\,\mathrm{d}s$.
\end{pf}

In \cite{nualartwhite}, it is proved that under some additional hypotheses, the mapping
$t\rightarrow B_t$ is differentiable from $(0,\infty)$ to $(\mathcal{S})^*$ and that its derivative,
denoted by $\dot B_t$, is a Hida distribution whose $S$-transform is $K_+\xi(t)$.

%s3.2 ###
\subsection{Bidimensional case}  \label{3.2}
Let $W=(W_t)_{t\in \R }$ be  a two-sided Brownian motion defined in the white noise probability space
$(\mathcal{S}'(\R ), \mathcal{B},  \bP)$. We can consider two independent standard Brownian motions as follows:
for $t\geq 0$, we set $W^{(1)}_t= W_t$ and $W^{(2)}_t = W_{-t}$.

In this section, we consider a bidimensional process $B=(B^{(1)}_t,B^{(2)}_t)_{t\geq 0}$,
where $B^{(1)}$ and $B^{(2)}$ are independent Volterra processes of the form
%e3.12 ###
\begin{equation}
B^{(i)}_t=\int_0^t K(t,s)\,\mathrm{d}W^{(i)}_s,  \qquad t\geq 0,  i=1,2.  \label{B1B2}
\end{equation}
For simplicity only, we work with the same kernel $K$ for the two components.

First, using exactly the same lines of reasoning as in the proof of Lemma~\ref{def-int},
we get the following result.
%l5
\begin{lem}\label{def-int-bis}
Let $B=(B^{(1)}_t,B^{(2)}_t)_{t\geq 0}$ be given as above,
with a kernel $K$ satisfying the conditions~$\mathrm{(H_1)}$ and $\mathrm{(H_2)}$. Assume, moreover, that $C$ defined by (\ref{C})
belongs to $L^2([0,T])$ for any $T>0$.
We then have the following results:
\begin{enumerate}[(1)]
\item[(1)] the
function $\xi\mapsto \int_0^T  (\int_0^u K_+\xi (-y)\,\mathrm{d}y ) K_+\xi (u) \, \mathrm{d}u$ is the $S$-transform of an
element of $(\mathcal{S})^*$, denoted by $\int_0^T B^{(1)}_u\diamond \dot{B}^{(2)}_u\,\mathrm{d}u$;
\item [(2)]the
function $\xi\mapsto \int_0^T K_+\xi (-u) K_+\xi (u) \,\mathrm{d}u$ is the $S$-transform of an
element of $(\mathcal{S})^*$, denoted by $\int_0^T \dot{B}^{(1)}_u\diamond \dot{B}^{(2)}_u\,\mathrm{d}u$.
\end{enumerate}
\end{lem}

We can now state the following result.
%p6
\begin{prop}\label{CVS*-bid}
Let $B=(B^{(1)}_t,B^{(2)}_t)_{t\geq 0}$ be given as above,
with a kernel $K$ satisfying the conditions
$\mathrm{(H_1)}$, $\mathrm{(H_2)}$ and $(\mathrm{H}_3^2)$.
The following convergences then hold:
\begin{eqnarray*}
\int_0^T B^{(1)}_u \frac{B^{(2)}_{u+\e}-B^{(2)}_u}{\e} \,\mathrm{d}u
&\displaystyle\mathop{\stackrel{  (\mathcal{S})^*}{\longrightarrow}}_{\e\to 0}  &
\int_0^T B^{(1)}_u\diamond \dot{B}^{(2)}_u \,\mathrm{d}u,  \\
\int_0^T  \biggl(\int_0^u \frac{B^{(1)}_{v+\e}-B^{(1)}_v}{\e}\,\mathrm{d}v \biggr)\frac{B^{(2)}_{u+\e}-B^{(2)}_u}{\e} \,\mathrm{d}u
&\displaystyle\mathop{\stackrel{  (\mathcal{S})^* }{\longrightarrow}}_{\e\to 0}&
\int_0^T B^{(1)}_u\diamond \dot{B}^{(2)}_u \,\mathrm{d}u,  \\
\int_0^T \frac{B^{(1)}_{u+\e}-B^{(1)}_u}{\e}\times\frac{B^{(2)}_{u+\e}-B^{(2)}_u}{\e} \,\mathrm{d}u
&\displaystyle\mathop{\stackrel{  (\mathcal{S})^*}{\longrightarrow}}_{\e\to 0}&
\int_0^T \dot{B}^{(1)}_u\diamond\dot{B}^{(2)}_u \,\mathrm{d}u.
\end{eqnarray*}
\end{prop}
\begin{pf}
Set
\[
\widetilde{G}_\e=
\int_0^T B^{(1)}_u \frac{B^{(2)}_{u+\e}-B^{(2)}_u}{\e} \,\mathrm{d}u=
\int_0^T B^{(1)}_u\diamond\frac{B^{(2)}_{u+\e}-B^{(2)}_u}{\e} \,\mathrm{d}u.
\]
From linearity and property (\ref{prod}) of the $S$-transform, we have
\[
S (\widetilde{G}_\e)(\xi)=\frac1\e  \int_0^T  S\bigl(B^{(1)}_u\bigr)(\xi )
S\bigl(B^{(2)}_{u+\e} - B^{(2)}_u\bigr) (\xi  )   \,\mathrm{d}u
\]
so that
\[
S (\widetilde{G}_\e)(\xi)=  \int_0^T  \biggl( \int_0^{u} K_+\xi (-y)\,\mathrm{d}y  \biggr)
\biggl( \frac1\e\int_u^{u+\e} K_+\xi (x)\,\mathrm{d}x  \biggr) \,\mathrm{d}u.
\]
Therefore, using (\ref{e11}) and (\ref{v12}), we can write
\begin{eqnarray*}
|S (\widetilde{G}_\e)(\xi)| &\leq&  M^2  |\xi|^2_{17/12}  \int_0^T  \biggl(\int_0^u C(t)\,\mathrm{d}t \biggr) D(u)\,\mathrm{d}u\\
&\leq& M^2  |\xi|^2_{17/12}  \int_0^T  \biggl(\int_0^u D(t)\,\mathrm{d}t \biggr)D(u)\,\mathrm{d}u\\
&=& \frac12 M^2  |\xi|^2_{17/12}  \biggl( \int_0^T D(u)\,\mathrm{d}u \biggr)^2\\
&\leq& \frac T2 M^2  |\xi|^2_{17/12}  \int_0^T D^2(u)\,\mathrm{d}u.
\end{eqnarray*}
Hence, by  the dominated convergence theorem, we get
%e3.13 ###
\begin{equation}\label{e12bis}
\lim_{\e\to 0} S (\widetilde{G}_\e)(\xi)=
\int_0^T  \biggl(\int_0^u K_+\xi (-y)\,\mathrm{d}y \biggr) K_+\xi (u)  \,\mathrm{d}u.
\end{equation}
The right-hand side of (\ref{e12bis})  is  the $S$-transform
of  $\int_0^T B^{(1)}_u \diamond \dot{B}^{(2)}_u\, \mathrm{d}u$, due to Lemma~\ref{def-int-bis}. Therefore, by
Theorem~\ref{th_kuo}, we obtain the desired result in point (1).

The proofs of the other two convergences follow exactly the same lines of reasoning and are therefore left to the reader.
\end{pf}

%s4 ###
\section{Fractional Brownian motion case}\label{sec4}

%s4.1 ###
\subsection{One-dimensional case}\label{sec31}

Consider a (one-dimensional) fractional Brownian motion (fBm)  $B=( B_t )_{t\geq 0}$ of Hurst index $H\in(0,1)$.
This means that $B$ is a zero mean Gaussian process with covariance function
\[
R_H(t,s)=E(B_tB_s) = \tfrac 12(t^{2H} + s^{2H}  - |t-s| ^{2H}).
\]
It is well known that $B$ is a Volterra process. More precisely (see \cite{DU}),
$B$ has the form (\ref{star}) with the kernel $K(t,s)=K_{H}(t,s)$ given by
\[
K_{H}(t,s)=c_{H} \biggl[  \biggl( \frac{t(t-s)}{s} \biggr) ^{H-1/2}
-\biggl(H-\frac{1}{2}\biggr)s^{1/2-H}\int_{s}^{t}u^{H-3/2}(u-s)^{H-1/2}\,\mathrm{d}u \biggr].
\]
Here, $c_H$ is  a constant depending only on $H$. Observe that
%e4.1 ###
\begin{equation}  \label{b1}
\frac{\partial K_{H}}{\partial t}(t,s)=c_{H}\biggl(H-\frac{1}{2}\biggr)(t-s)^{H-\sfrac{3}{2}} \biggl( \frac{s}{t} \biggr) ^{\sfrac{1}{2}-H}\qquad \mbox{for $t>s> 0$.}
\end{equation}

Denote by $\mathscr{E}$ the set of all $\R$-valued
step functions defined on $[0,\infty)$. Consider the Hilbert space $\HH$ obtained by closing $\mathscr{E}$
with respect to the inner product
\[
\big\langle \mathbf{1}_{[0,u]},\mathbf{1}_{[0,v]}\big\rangle_\HH=E(B_uB_v).
\]
The mapping $\mathbf{1}_{[0,t]}\mapsto B_t$ can be extended to an isometry $\varphi\mapsto B(\varphi)$ between $\HH$
and the Gaussian space $\mathcal{H}_1$ associated with $B$.
Also, write $\HH^{\otimes k}$ to indicate the $k$th tensor product of $\HH$.
When $H>1/2$,  the inner product in the space $\HH$ can be written as follows, for any
$\varphi$, $\psi\in \mathscr{E}$:
\[
\langle \phi, \psi \rangle _{\HH}=H(2H-1)\int_0^\infty\!\! \int_0^\infty \phi(s) \psi(s') |s-s'|^{2H-2}\,\mathrm{d}s\,\mathrm{d}s'.
\]
By approximation, this extends immediately to any $\varphi$, $\psi\in {\mathcal{S}}(\R)\cup\mathscr{E}$.

We will make use of  the multiple integrals with respect to $B$
(we refer to \cite{nualart} for a detailed account on the properties of these integrals).
For every $k\geq 1$, let $\mathcal{H}_{k}$ be the $k$th Wiener chaos of $B$,
that is, the closed linear subspace of $L^{2}(\Omega)$
generated by the random variables $\{h_{k}(B(\varphi)), \varphi\in \mathfrak{H},  \|
\varphi \| _{\mathfrak{H}}=1\}$, where $h_{k}$ is the $k$th Hermite polynomial (\ref{herm-pol}).
For
any $k\geq 1$, the mapping $I_{k}(\varphi^{\otimes k})=h_{k}(B(\varphi))$ provides a
linear isometry between the \textit{symmetric} tensor product $\mathfrak{H}^{\odot k}$
(equipped with the modified norm $\sqrt{k!} \| \cdot  \| _{\mathfrak{H}^{\otimes k}}$) and the $k$th Wiener chaos~$\mathcal{H}_{k}$.\looseness=1

Following \cite{NNT}, let us now introduce the Hermite random variable $Z^{(k)}_T$ mentioned in (\ref{cv>}).
Fix $T>0$ and let $k\geq 1$ be an integer.
The family $(\varphi_\e)_{\e>0}$, defined by
%e4.2 ###
\begin{equation}\label{phieps}
\varphi _{\e}=\e^{-k}\int_0^T \mathbf{1}_{[u,u+\e]}^{\otimes k}\,\mathrm{d}u,
\end{equation}
satisfies
%e4.3 ###
\begin{eqnarray}\label{variance}
\hspace*{-20pt}\lim_{\e,\eta\rightarrow 0 } \langle \varphi_\e,
\varphi _{\eta} \rangle _{\mathfrak{H}^{\otimes k}}
=
H^k(2H-1)^{k}
\int_{[0,T]^2}|s-s'|^{(2H-2)k}\,\mathrm{d}s\,\mathrm{d}s'
=c_{k,H} T^{(2H-2)k+2}
\end{eqnarray}
with $c_{k,H}=\frac{H^k(2H-1)^{k}}{(Hk-k+1)(2Hk-2k+1) }$.
This implies that  $\varphi_\e$ converges, as $\e$ tends to zero, to an element of
$\mathfrak{H}^{\otimes k}$. The limit, denoted by $\pi^k_{\mathbf{1}_{[0,T]}}$,
can be characterized as follows.
For any $\xi_i\in {\mathcal{S}}(\R)$, $i=1,\ldots,k$, we have
\begin{eqnarray*}
&&\big\langle\pi^k_{\mathbf{1}_{[0,T]}},\xi_1\otimes \cdots\otimes \xi_k\big\rangle_{\mathfrak{H}^{\otimes k}}\\
&&\quad=\lim_{\e\to 0} \langle \varphi_\e,\xi_1\otimes \cdots\otimes \xi_k\rangle_{\mathfrak{H}^{\otimes k}}\\
&&\quad=\lim_{\e\to 0} \e^{-k}\int_0^T \mathrm{d}u \,\prod_{i=1}^k \big\langle \mathbf{1}_{[u,u+\e]},\xi_i\big\rangle_{\mathfrak{H}}\\
&&\quad=\lim_{\e\to 0} \e^{-k}H^k(2H-1)^k\int_0^T \mathrm{d}u\, \prod_{i=1}^k\int_{u}^{u+\e} \mathrm{d}s\, \int_0^T
\mathrm{d}r\,
|s-r|^{2H-2}\xi_i(r)\\
&&\quad=
H^k(2H-1)^k \int_0^T \mathrm{d}u\, \prod_{i=1}^k \int_0^T \mathrm{d}r\,|u-r|^{2H-2}\xi_i(r).
\end{eqnarray*}
We define the $k$th Hermite random variable
by $Z^{(k)}_T=
I_{k}(\pi^k_{\mathbf{1}_{[0,T]}})$.
Note that, by using the isometry formula for multiple integrals
and since $G_\e=I_k(\varphi_\e)$,
the convergence (\ref{cv>}) is just a corollary of our construction of $Z^{(k)}_T$.
Moreover, by (\ref{variance}), we have
\[
E \bigl[\bigl(Z_T^{(k)}\bigr)^2 \bigr]=c_{k,H}\times t^{(2H-2)k+2}.
\]

We will need the following preliminary result.
\eject
%l7
\begin{lem}\label{lmlm}
\begin{enumerate}[(1)]
\item [(1)]The fBm $B$ verifies the assumptions $\mathrm{(H_1)}$, $\mathrm{(H_2)}$ and $(\mathrm{H}_3^{k})$ if and only if $H\in (\frac12-\frac1k,1 )$.
\item [(2)]If $H\in (\frac12-\frac1k,1 )$, then $\int_0^T \dot{B}_u^{\diamond k}\,\mathrm{d}u$ is a well-defined element of $(\mathcal{S})^*$ (in the sense
of Lemma~\ref{def-int}).
\item [(3)]If we assume that $H>\frac12$, then $\int_0^T \dot{B}_u^{\diamond k}\,\mathrm{d}u$ belongs to
$L^2(\Omega)$ if and only if $H>1-\frac1{2k}$.
\end{enumerate}
\end{lem}
\begin{pf}
(1) Since
%e4.4 ###
\begin{equation}\label{k'}
k'(t)=k'_H(t)=\bigl(H+ \tfrac 12 \bigr)c_1(H)t^{H-\sfrac 12 }
\end{equation}
and
%e4.5 ###
\begin{equation}\label{terme2}
\int_0^t  \bigg|\frac{\partial K_H}{\partial t}(t,s) \bigg|(t-s)\,\mathrm{d}s=
\bigg|\int_0^t \frac{\partial K_H}{\partial t}(t,s)(t-s)\,\mathrm{d}s \bigg|=
c_2(H)t^{H+ \sfrac 12}
\end{equation}
for some constants $c_1(H)$ and $c_2(H)$,
we immediately see that assumptions $\mathrm{(H_1)}$ and $\mathrm{(H_2)}$ are satisfied for all $H\in(0,1)$. It therefore remains
to focus on assumption $(\mathrm{H}_3^{k})$.
For all $H\in(0,1)$, we have
%e4.6 ###
\begin{equation}\label{$1}
\sup_{0<\e\leq \e_0}\frac1\e\int_{t}^{t+\e} s^{H-1/2}\,\mathrm{d}s\leq t^{H-\sfrac 12} \vee
(t+\e_0)^{H-\sfrac 12}
\end{equation}
and
%e4.7 ###
\begin{equation}\label{$3}
\sup_{0<\e\leq \e_0}\frac1\e\int_{t}^{t+\e} s^{H+1/2}\,\mathrm{d}s\leq (t+\e_0)^{H+1/2}.
\end{equation}
Consequently, since $\int_0^T t^{kH-k/2}\,\mathrm{d}t$ is finite when $H>\frac12-\frac1k$, we deduce
from (\ref{k'})--(\ref{$3})
that
$(\mathrm{H}_{3}^{k})$ holds in this case. Now, assume that $H\leq \frac12-\frac1k$.
Using the fact that $D(t) \ge C(t)$, we obtain
\[
\int_0^T D^k(t)\,\mathrm{d}t \ge  \int_0^T C^k(t)\,\mathrm{d}t
= \biggl(H+\frac 12\biggr)^kc_1(H)^k \int_0^T t^{kH-\sfrac k2}\,\mathrm{d}t=\infty.
\]
Therefore, in this case, assumption $(\mathrm{H}_{3}^{k})$ is not verified.

(2) This fact can be proven immediately: simply combine  the previous point with Lemma~\ref{def-int}.

(3) By definition of $\int_0^T \dot{B}_u^{\diamond k}\, \mathrm{d}u$ (see Lemma~\ref{def-int}),
it is equivalent to show that the distribution $\tau^{k}_{\1_{[0,T]}}$, defined via the identity
$\int_0^t \dot{B}_s^{\diamond k}\,\mathrm{d}s= I_k(\tau^{k}_{\1_{[0,t]}})$,
can be represented as a
function belonging to $L^2([0,T]^k)$.
We can write
\begin{eqnarray*}
\big\langle \tau^{k}_{\1_{[0,T]}}, \xi_1\otimes\cdots \otimes \xi_k \big\rangle & =&
\int_0^T K_+\xi_1(s)\cdots K_+\xi_k(s)\,\mathrm{d}s\\
& =& \int_0^T \mathrm{d}s\,\prod_{i=1}^{k} \int_0^s \frac{\partial K_H}{\partial s}(s,r)\xi_i(r)\,\mathrm{d}r
\end{eqnarray*}
for any  $\xi_1,\ldots, \xi_k\in \mathcal{S}(\R)$.
Observe that $K_+\xi(s)=\int_0^s \frac{\partial K_H}{\partial s}(s,r)\xi(r)\,\mathrm{d}r$
because $K_H(s,s)=0$ for $H>1/2$.
Using Fubini's theorem, we deduce that the distribution $\tau^{k}_{\1_{[0,T]}}$ can be represented as the
function
\begin{eqnarray*}
 \tau _{\mathbf{1}_{[0,T]}}^{k}(x_{1},\ldots,x_{k})=\mathbf{1}_{[0,T]^{k}}(x_{1},\ldots,x_{k})
\int_{\max (x_{1},\ldots,x_{k})}^{T}\frac{\partial K_H}{\partial s}(s,x_{1})\cdots \frac{\partial K_H}{\partial s}
(s,x_{k})\,\mathrm{d}s.
\end{eqnarray*}
We then obtain
\begin{eqnarray*}
&&\big\| \tau _{\mathbf{1}_{[0,T]}}^{k}\big \| _{L^2([0,T]^k)}^{2}
\\
&&\quad=\int_{[0,T]^{k}}\int_{\max (x_{1},\ldots,x_{k})}^{T}\int_{\max (x_{1},\ldots,x_{k})}^{T}\frac{\partial K_H}
{\partial s}(s,x_{1})\cdots \frac{\partial K_H}{\partial s}(s,x_{k}) \\
&&\qquad\hspace*{130pt}{}\times \frac{\partial K_H}{\partial s}(r,x_{1})\cdots \frac{\partial K_H}{\partial s}(r,x_{k})
\,\mathrm{d}s\,\mathrm{d}r\,\mathrm{d}x_{1}\cdots \,\mathrm{d}x_{k} \\
&&\quad=\int_{[0,T]^{2}} \biggl( \int_{0}^{r\wedge s}\frac{\partial K_H}{\partial s}(s,x)\frac{\partial K_H}{\partial s}(r,x)\,\mathrm{d}x \biggr)^k
\,\mathrm{d}r\,\mathrm{d}s.
\end{eqnarray*}
Using the equality  (\ref{b1}) and the same computations as in \cite{nualart}, page 278, we obtain, for $s<r$,
%e4.8 ###
\begin{equation}\label{etunetunetunzero}
\int_{0}^{s}\frac{\partial K_H}{\partial s}(s,x)\frac{\partial K_H}{\partial r}   (r,x)\,\mathrm{d}x=
H(2H-1) (r-s)^{2H-2}.
\end{equation}
Therefore,
\[
\big\|\tau^k_{\mathbf{1}_{[0,T]}}\big \|^2_{L^2([0,T]^k)}=
\bigl(H(2H-1)\bigr)^k  \int_0^T\int_0^T  |r-s|^{2Hk-2k} \,\mathrm{d}r\,\mathrm{d}s.
\]
We immediately check that  $\|\tau^k_{\mathbf{1}_{[0,T]}}  \|^2_{L^2([0,T]^k)}<\infty$ if and only if
$2Hk-2k>-1$, that is,
$H>1-\frac{1}{2k}$. Thus, in this case, the Hida distribution
$\int_0^T \dot{B}_s^{\diamond k}\,\mathrm{d}s$ is a square-integrable random variable with
\begin{eqnarray*}
E \biggl[ \biggl(\int_0^T \dot{B}_s^{\diamond k}\,\mathrm{d}s  \biggr)^2 \biggr] =
\big\|\tau^k_{\mathbf{1}_{[0,T]}}\big\|^2_{L^2([0,T]^k)}
=c_{k,H}\times T^{2Hk-2k+2}.
\end{eqnarray*}
\upqed
\end{pf}

%r8
\begin{rem}
According to our result, the two distributions
$\tau^k_{\mathbf{1}_{[0,T]}}$ and $\pi^k_{\mathbf{1}_{[0,T]}}$ should coincide when
$H>1/2$. We can check this fact by means of elementary arguments. Let $\xi_i\in {\mathcal{S}}(\R)$, $i=1,\ldots,k$. From (\ref{e7}), we deduce that
\[
\big\langle \mathbf{1}_{[u,u+\e]},\xi_i\big\rangle_{\mathfrak{H}}=\int_u^{u+\e}K_+\xi_i(s) \,\mathrm{d}s
\]
and then
\[
\lim_{\e\to 0} \frac{1}{\e}\big\langle \mathbf{1}_{[u,u+\e]},\xi_i\big\rangle_{\mathfrak{H}}= K_+\xi_i(u).
\]
Using (\ref{v12}) with $k=1$ for each $\xi_i$ and applying the dominated convergence theorem,
since the fractional Brownian motion satisfies the assumption $(\mathrm{H}_3^k)$ when
$H\in (\frac12-\frac1k,1 )$, we get, for $\varphi_\e$ defined in (\ref{phieps}),
\[
\lim_{\e\to 0} \langle \varphi_\e,\xi_1\otimes \cdots\otimes \xi_k\rangle_{\mathfrak{H}^{\otimes k}}
= \int_0^T K_+\xi_1(u)\cdots K_+\xi_k(u)\,\mathrm{d}u,
\]
which yields $\tau^k_{\mathbf{1}_{[0,T]}}=\pi^k_{\mathbf{1}_{[0,T]}}$.
\end{rem}

We can now state the main result of this section.
%p9
\begin{prop}\label{thm-fbm}
Let $k\geq 2$ be an integer. If $H>\frac 12- \frac 1k$ (note that this condition is immaterial for $k=2$), the  random variable
\[
G_\e=\int_0^T  \biggl(\frac{B_{u+\e}-B_u}{\e} \biggr)^{\diamond k}\,\mathrm{d}u
=\e^{-k(1-H)}\int_0^T h_k
\biggl(\frac{B_{u+\e}-B_u}{\e^{H}}\biggr)\,\mathrm{d}u
\]
converges in $(\mathcal{S}^*)$, as $\e\to 0$, to the Hida distribution
$\int_0^T \dot{B}_u^{\diamond k}\,\mathrm{d}u$.
Moreover, $G_\e$ converges in $L^2(\Omega)$ if and only if $H>1-\frac{1}{2k}$.
In this case, the limit is $\int_0^T \dot{B}_u^{\diamond k}\,\mathrm{d}u=Z_T^{(k)}$.
\end{prop}

\begin{pf}
The first point follows directly from Proposition~\ref{CVS*} and Lemma~\ref{lmlm} (point 1).
On the other hand, we already know (see (\ref{cv>})) that $G_\e$ converges in $L^2(\Omega)$
to $Z_T^{(k)}$ when $H>1-\frac1{2k}$. This implies that when $H>1-\frac1{2k}$,
$\int_0^T \dot{B}_s^{\diamond k}\,\mathrm{d}s$ must be a square-integrable random variable equal to $Z_T^{(k)}$.
Assume, now, that $H\leq 1-\frac1{2k}$. From the proof of (\ref{cv<}) and (\ref{cv=})
below, it follows that $E(G^2_\e)$ tends to $+\infty$ as $\e$ tends to zero,
so $G_\e$ does not converge in $L^2(\Omega)$.
\end{pf}

%s4.2 ###
\subsection{Bidimensional case}\label{sec32}
Let $B^{(1)}$ and $B^{(2)}$ denote two independent fractional Brownian motions with (the same) Hurst index $H\in(0,1)$,
defined by the stochastic integral representation
(\ref{B1B2}), as in Section~\ref{3.2}.

By combining Lemma~\ref{lmlm} (point 1 with $k=2$) and Lemma~\ref{def-int-bis},
we have the following preliminary result.
%l10
\begin{lem}\label{lmlmlm}
For all $H\in(0,1)$, the Hida distributions
$\int_0^T B^{(1)}_u \diamond \dot{B}^{(2)}_u\,\mathrm{d}u$
and $\int_0^T \dot{B}^{(1)}_u \diamond \dot{B}^{(2)}_u\,\mathrm{d}u$ are well-defined elements
of $(\mathcal{S})^*$ (in the sense
of Lemma~\ref{def-int-bis}).
\end{lem}

We can now state  the following result.
%p11
\begin{prop}\label{thm-fbm2}
\begin{enumerate}[(1)]
\item[(1)] For all $H\in(0,1)$, $\widetilde{G}_\e$ defined by (\ref{tilde})
converges in $(\mathcal{S}^*)$, as $\e\to 0$, to the Hida distribution
$\int_0^T B^{(1)}_u\diamond \dot{B}^{(2)}_u \,\mathrm{d}u$.
Moreover, $\widetilde{G}_\e$ converges in $L^2(\Omega)$ if and only if $H\geq 1/2$.
\item[(2)] For all $H\in(0,1)$, $\breve{G}_\e$ defined by (\ref{breve})
converges in $(\mathcal{S}^*)$, as $\e\to 0$, to the Hida distribution
$\int_0^T B^{(1)}_u\diamond \dot{B}^{(2)}_u \,\mathrm{d}u$.
Moreover, $\breve{G}_\e$ converges in $L^2(\Omega)$ if and only if $H> 1/4$.
\item[(3)] For all $H\in(0,1)$, $\hat{G}_\e$
defined by (\ref{hat})
converges in $(\mathcal{S}^*)$, as $\e\to 0$, to the Hida distribution
$\int_0^T \dot{B}^{(1)}_u\diamond \dot{B}^{(2)}_u \,\mathrm{d}u$.
Moreover, $\widehat{G}_\e$ converges in $L^2(\Omega)$ if and only if $H>3/4$.
\end{enumerate}
\end{prop}

\begin{pf}
(1)
The first point follows directly from Proposition~\ref{CVS*-bid} and Lemma~\ref{lmlm} (point 1 with $k=2$).
Assume that $H<1/2$.
From the proof of Theorem~\ref{cvgausslevyarea} below, it follows that $E(\widetilde{G}^2_\e)\to \infty$ as
$\e$ tends to zero, so $\widetilde{G}_\e$ does not converge in $L^2(\Omega)$.
Assume that $H=1/2$. By a classical result of Russo and Vallois (see, for example, the survey \cite{RVLN})
and since we are, in this case, in a martingale setting, we have that $\widetilde{G}_\e$
converges in $L^2(\Omega)$ to the It\^o integral $\int_0^T B^{(1)}_u \,\mathrm{d}B^{(2)}_u$.
Finally, assume that $H> 1/2$.
For $\e,\eta>0$, we have
\[
E( \widetilde{G}_\e \widetilde{G}_\eta)
= \frac 1{\e \eta}  \int_{[0,T]^2} \rho_{\e,\eta}(u- u'  ) R_H(u,u' ) \,\mathrm{d}u\,\mathrm{d}u',
\]
where
%e4.9 ###
\begin{equation}\label{rhoeps}
\rho_{\e,\eta}(x )= \tfrac 12 [  |x+\e|^{2H} + |x-\eta|^{2H} - |x|^{2H}
-|x+\e-\eta |^{2H} ].
\end{equation}
Note that as $\e$ and $\eta$ tend to zero, the quantity
$(\e\eta)^{-1} \rho_{\e,\eta}(u- u'  )$ converges pointwise to (and is bounded by)  $H(2H-1) |u-u' |^{2H-2}$.
Then, by the dominated convergence theorem, it follows that  $E( \widetilde{G}_\e \widetilde{G}_\eta)$ converges to
\[
H(2H-1) \int_{[0,T]^2} |u-u' |^{2H-2}R_H(u,u' ) \,\mathrm{d}u\,\mathrm{d}u'
\]
as $\e,\eta\to 0$, with $\int_{[0,T]^2} |u-u' |^{2H-2}|R_H(u,u' )| \,\mathrm{d}u\,\mathrm{d}u' <\infty$, since $H> 1/2$.
Hence, $\widetilde{G}_\e$ converges in $L^2(\Omega)$.

(2)
The first point follows directly from Proposition~\ref{CVS*-bid} and Lemma~\ref{lmlm} (point 1 with $k=2$).
Assume that $H\leq 1/4$.
From the proof of Theorem~\ref{cvgausslevyareabis} below, it follows that $E(\breve{G}^2_\e)\to \infty$ as
$\e$ tends to zero, so $\breve{G}_\e$ does not converge in $L^2(\Omega)$.
Assume that $H>1/4$.
For $\e,\eta>0$, we have
\[
E( \breve{G}_\e \breve{G}_\eta)
= \frac 1{\e^2 \eta^2}  \int_{[0,T]^2}\,\mathrm{d}u\,\mathrm{d}u'  \rho_{\e,\eta}(u- u'  )
\int_0^u\mathrm{d}s\int_0^{u' }\mathrm{d}s'\, \rho_{\e,\eta}(s-s')
\]
with $\rho_{\e,\eta}$ given by (\ref{rhoeps}).
Note that, as $\e$ and $\eta$ tend to zero, the quantity
$(\e\eta)^{-1} \rho_{\e,\eta}(u- u'  )$ converges pointwise to  $H(2H-1) |u-u' |^{2H-2}$,
whereas $ (\e\eta)^{-1}\int_0^u\mathrm{d}s\int_0^{u' }\mathrm{d}s'\,\rho_{\e,\eta}(s-s')$ converges pointwise
to $R_H(u,u' )$.
It then follows that  $E(\breve{G}_\e \breve{G}_\eta)$ converges to
\[
-\frac{H}2(2H-1) \int_{[0,T]^2}|u-u' |^{4H-2}\,\mathrm{d}u\,\mathrm{d}u' +H\int_0^T u^{2H}\bigl(u^{2H-1}+(T-u)^{2H-1}\bigr)\,\mathrm{d}u
\]
as $\e,\eta\to 0$ and each integral is finite since $H> 1/4$.
Hence, $\breve{G}_\e$ converges in $L^2(\Omega)$.

(3) Once again, the first point follows
from Proposition~\ref{CVS*-bid} and Lemma~\ref{lmlm} (point 1 with $k=2$).
Assume that $H\leq 3/4$. From the proof of Theorem~\ref{cvgausscovariation} below,
it follows that $E(\widehat{G}^2_\e)\to \infty$ as
$\e$ tends to zero, so $\widehat{G}_\e$ does not converge in $L^2(\Omega)$.
Assume, now, that $H>3/4$.
For $\e,\eta>0$, we have
\[
E( \widehat{G}_\e \widehat{G}_\eta)
= \frac 1{\e^2 \eta^2}  \int_{[0,T]^2} \rho_{\e,\eta}(u- u'  )^2 \,\mathrm{d}u\,\mathrm{d}u'
\]
with $\rho_{\e,\eta}$ given by (\ref{rhoeps}).
Since  the quantity
$(\e\eta)^{-1} \rho_{\e,\eta}(u- u'  )$ converges pointwise to (and is bounded by) $H(2H-1) |u-u' |^{2H-2}$,
we have, by the dominated convergence theorem, that
$E( \widehat{G}_\e \widehat{G}_\eta)$ converges to
\[
H^2(2H-1)^2 \int_{[0,T]^2} |u-u' |^{4H-4}\,\mathrm{d}u\,\mathrm{d}u'
\]
as $\e,\eta\to 0$,
with $\int_{[0,T]^2} |u-u' |^{4H-4}\,\mathrm{d}u\,\mathrm{d}u' <\infty$, since $H>3/4$.
Hence, $\widehat{G}_\e$ converges in $L^2(\Omega)$.
\end{pf}

%s5 ###
\section[Proof of  the convergences (1.3) and (1.4)]{Proof of
the convergences (\protect\ref{cv<}) and (\protect\ref{cv=})}\label{sec5}

In  this section, we provide a new proof of these convergences by means of a
recent criterion for the weak convergence of sequences of multiple stochastic integrals
established in \cite{NP} and \cite{PT}. We refer to  \cite{MR1} for a proof in the case of more general Gaussian processes, using different kind of tools.

Let us first recall the aforementioned criterion.
We continue to use the notation introduced in Section~\ref{sec31}.
Also, let $\{e_{i}, i\geq 1\}$ denote a complete orthonormal system in $\mathfrak{H}$.
Given $f\in \mathfrak{H}^{\odot k}$ and $g\in \mathfrak{H}^{\odot l}$, for every $r=0,\ldots,k\wedge l$, the \textit{contraction} of $f$ and $g$ of order $r$
is the element of $\mathfrak{H}^{\otimes (k+l-2r)}$ defined by
\[
f\otimes _{r}g=\sum_{i_{1},\ldots,i_{r}=1}^{\infty }\langle
f,e_{i_{1}}\otimes \cdots \otimes e_{i_{r}}\rangle _{\mathfrak{H}^{\otimes
r}}\otimes \langle g,e_{i_{1}}\otimes \cdots \otimes e_{i_{r}}\rangle _{\mathfrak{H}^{\otimes r}}.
\]
(Note that $f\otimes_{0}g=f\otimes g$ equals the tensor product of $f$ and
$g$ while, for $k=l$, $f\otimes _{k}g=\langle f,g\rangle _{\mathfrak{H}^{\otimes k}}$.) Fix $k\geq 2$ and let $(F_\e)_{\e>0}$ be a family of the form
$F_\e=I_k(\phi_\e)$ for some $\phi_\e\in\HH^{\odot k}$.
Assume that\vspace*{-2pt} the variance of $F_\e$ converges as $\e\to 0$ (to $\sigma^2$, say).
The criterion of Nualart and Peccati \cite{NP} asserts that
$F_\e\stackrel{\mathrm{Law}}{\longrightarrow}N\sim\mathscr{N}(0,\sigma^2)$
if and only if
$  \| \phi_\e \otimes_r  \phi_\e\|_{\HH^{\otimes (2k-2r)} } \to 0$
for any $r=1,\ldots, k-1$.
In this case,  due to the result proved by Peccati and Tudor \cite{PT},
we automatically have  that
\[
(B_{t_1},\ldots,B_{t_k},F_\e)\stackrel{\mathrm{Law}}{\longrightarrow}(B_{t_1},\ldots,B_{t_k},N)
\]
for all $t_k> \cdots> t_1>0$, with $N\sim\mathscr{N}(0,\sigma^2)$ \textit{independent} of $B$.

For $x\in\R$, set
%e5.1 ###
\begin{equation} \label{rho}
\rho(x)=  \tfrac 12  (  |x+1|^{2H} + |x-1|^{2H} - 2|x|^{2H} ),
\end{equation}
and note that $\rho(u-v)=E [(B_{u+1}-B_u)(B_{v+1}-B_v) ]$ for all $u,v\geq 0$
and that $\int_\R |\rho(x)|^k \,\mathrm{d}x$ is finite if and only if $H<1-\frac1{2k}$ (since
$\rho(x)\sim H(2H-1)|x|^{2H-2}$ as $|x|\to\infty$).

We now proceed with the proof of  (\ref{cv<}). The proof  of  (\ref{cv=}) would follow similar arguments.

\begin{pf*}{Proof of (\ref{cv<})}
Because $\e^{k(1-H)-\sfrac12}G_\e$ can be expressed as  a $k$th multiple Wiener integral, we can use the criterion of Nualart and Peccati.
By the scaling property of the fBm, it is actually equivalent to considering the family of random variables
$(F_\e)_{\e>0}$, where
\[
F_\e = \sqrt{\e}  \int_0^{T/\e}  h_k(B_{u+1} - B_u) \,\mathrm{d}u.
\]

\textit{Step \textup{1}. Convergence of the variance}.
We can write
\begin{eqnarray*}
E(F^2_{\e})&=&    \e k!
\int_0^{T/\e}  \mathrm{d}u \int_0^{T/\e} \mathrm{d}s\,  \rho(u-s)^k
\\
 &=&  \e k! \int_{-T/\e}^{T/\e} \rho(x)^{k}(T/\e - |x|)\,\mathrm{d}x,
\end{eqnarray*}
where the function $\rho$ is defined in (\ref{rho}).  Therefore, by the dominated convergence theorem,
\[
\lim_{\e \downarrow 0} E(F^2_{\e})
=T  k!  \int_{\R} \rho(x)^{k}\,\mathrm{d}x.
\]

\textit{Step \textup{2}. Convergence of the contractions}.
Observe that
the random variable $h_{k }(B_{u+1} -B_u)$ coincides with the multiple stochastic integral
$I_k( \mathbf{1}_{[u,u+1]} ^{\otimes  k  })$.  Therefore,
$
F_\e = I_k(\phi_\e),
$
where
$   \phi_\e =   \sqrt{\e}    \int_0^{T/\e}   \mathbf{1}_{[u,u+1]} ^{\otimes  k  }\,\mathrm{d}u.$
Let $r\in\{1,\ldots,k-1\}$. We have
\[
\phi_\e \otimes_r  \phi_\e=  \e   \int_0^{T/\e} \int_0^{T/\e}   \bigl( \mathbf{1}_{[u,u+1]} ^{\otimes  (k-r)  }
\otimes  \mathbf{1}_{[s,s+1]} ^{\otimes ( k-r)  }  \bigr)  \rho(u-s)^r \,\mathrm{d}u\,\mathrm{d}s.
\]
As a consequence, $\| \phi_\e \otimes_r  \phi_\e\|_{\HH^{\otimes (2k-2r)} }^2$ equals
\[
\e^2   \int_{[0,{T/\e}]^4}  \rho(u-s)^r \rho(u' -s')^r
\rho(u-u' )^{k-r}    \rho(s-s')^{k-r}  \,\mathrm{d}s\,\mathrm{d}s'\,\mathrm{d}u\,\mathrm{d}u'.
\]
Making the changes of variables $x=u-s $, $y=u' -s'$ and $z=u-u' $, we obtain that
$    \| \phi_\e \otimes_r  \phi_\e\|_{\HH^{\otimes (2k-2r)} }^2$ is less than
\[
A_\e= \e    \int_{D_\e}  |\rho(x)|^r |\rho(y)|^r
|\rho(z)|^{k-r}    |\rho(y+z-x)|^{k-r}  \,\mathrm{d}x\,\mathrm{d}y\,\mathrm{d}z,
\]
where $D_\e= [-T/\e,{T/\e}]^3$. Consider the decomposition
\begin{eqnarray*}
A_\e&=& \e   \int_{D_\e \cap \{|x|\vee |y| \vee |z| \leq K\}}  |\rho(x)|^r |\rho(y)|^r
|\rho(z)|^{k-r}    |\rho(y+z-x)|^{k-r}  \,\mathrm{d}x\,\mathrm{d}y\,\mathrm{d}z  \\
&&{}   +\e   \int_{D_\e \cap \{|x|\vee |y| \vee |z| > K\}}  |\rho(x)|^r |\rho(y)|^r
|\rho(z)|^{k-r}    |\rho(y+z-x)|^{k-r}  \,\mathrm{d}x\,\mathrm{d}y\,\mathrm{d}z \\
&=& B_{\e,K} + C_{\e,K}.
\end{eqnarray*}
Clearly, for any fixed $K>0$, the term $B_{\e,K}$ tends to zero because $\rho$ is a bounded function.
On the other hand, we have
\[
D_\e \cap \{|x|\vee |y| \vee |z| > K\}\subset
D_{\e,K,x}\cup D_{\e,K,y}\cup D_{\e,K,z},
\]
where $D_{\e,K,x}=\{|x|>K\}\cap \{|y|\leq T/\e\}\cap \{|z|\leq T/\e\}$ ($D_{\e,K,y}$ and $D_{\e,K,z}$ being defined similarly). Set
\[
C_{\e,K,x}
=\e   \int_{D_{\e,K,x}}  |\rho(x)|^r |\rho(y)|^r
|\rho(z)|^{k-r}    |\rho(y+z-x)|^{k-r}  \,\mathrm{d}x\,\mathrm{d}y\,\mathrm{d}z.
\]
By H\"older's inequality, we have
\begin{eqnarray*}
C_{\e,K,x} &\leq& \e    \biggl(  \int_{D_{\e,K,x}}  |\rho(x)|^k|\rho(y)|^k \,\mathrm{d}x\,\mathrm{d}y\,\mathrm{d}z  \biggr)^{\sfrac rk}\\
&&{}\times
\biggl(  \int_{D_{\e,K,x}}  |\rho(z)|^k|\rho(y+z-x)|^k   \,\mathrm{d}x\,\mathrm{d}y\,\mathrm{d}z \biggr)^{1-\sfrac{r}k}\\
&\leq&
2T \biggl(\int_\R |\rho(t)|^k \,\mathrm{d}t \biggr)^{2-\sfrac rk} \biggl(\int_{|x|>K}|\rho(x)|^k \,\mathrm{d}x \biggr)^{\sfrac rk}
\mathop{\longrightarrow}_{K\to\infty} 0.
\end{eqnarray*}
Similarly, we prove that $C_{\e,K,y}\to 0$
and $C_{\e,K,z}\to 0$ as $K\to\infty$.
Finally, it suffices to choose $K$ large enough in order to get the desired result,
that is, $\| \phi_\e \otimes_r  \phi_\e\|_{\HH^{\otimes (2k-2r)} }\to 0$ as $\e\to 0$.

\textit{Step \textup{3}. Proof of the first point.}
By step 1, the family
\[
\bigl((B_t)_{t\in[0,T]},\e^{\sfrac12-2H}G_\e \bigr)
\]
is tight in $C([0,T])\times \mathbb{R}$.
By step 2,
we also have the convergence of the finite-dimensional distributions,
as a by-product of the criteria of Nualart and Peccati \cite{NP}
and Peccati and Tudor \cite{PT} (see the preliminaries at the beginning of this section).
Hence, the proof of the first point is complete.
\end{pf*}

%s6 ###
\section{Convergences in law for some functionals related to the L\'evy area of the fractional Brownian
motion}\label{sec6}

Let $B^{(1)}$ and $B^{(2)}$ denote two independent fractional Brownian motions
with Hurst index $H\in(0,1)$.
Recall the definition (\ref{tilde}) of $\widetilde{G}_\e$:
\[
\widetilde{G}_\e=\int_0^T B^{(1)}_u \frac{B^{(2)}_{u+\e}-B^{(2)}_u}{\e}\, \mathrm{d}u.
\]

%t12
\begin{theo}\label{cvgausslevyarea}
Convergence in law (\ref{star1}) holds.
\end{theo}

\begin{pf}
We fix $H<1/2$. The proof is divided into several steps.

\textit{Step \textup{1}. Computing the variance of  $\varepsilon ^{\sfrac{1}{2}-H}\widetilde{G}_\e$}.

By using the scaling properties
of the fBm,  first observe that $\varepsilon ^{\sfrac{1}{2}-H}\widetilde{G}_\e$ has the same law
as
%e6.1 ###
\begin{equation}  \label{gg}
\widetilde{F}_{\varepsilon }=\e^{1/2+H} \int_{0}^{T/\varepsilon
}B_{u}^{(1)} \bigl( B_{u+1}^{(2)}-B_{u}^{(2)} \bigr)\, \mathrm{d}u.
\end{equation}
For $\rho(x)=\frac12 (|x+1|^{2H}+|x-1|^{2H}-2|x|^{2H} )$, we have
\begin{eqnarray*}
E(\widetilde{F}_{\varepsilon }^{2})
&=&\varepsilon^{1+2H} \int_{0}^{T/\varepsilon}\mathrm{d}u\int_{0}^{T/\varepsilon }\mathrm{d}s\, R_{H}(u,s)\rho (u-s)
\\
&=&\alpha_{\varepsilon }-\beta_{\varepsilon },
\end{eqnarray*}
where
\begin{eqnarray*}
\alpha_{\varepsilon }
&=& \varepsilon^{1+2H} \int_{0}^{T/\varepsilon}\mathrm{d}u\, u^{2H}\int_{0}^{T/\e}\mathrm{d}s\, \rho (u-s),
\\
\beta_{\varepsilon }
&=&\varepsilon^{1+2H} \int_{0}^{T/\varepsilon}\mathrm{d}u\,\int_{0}^{u}\mathrm{d}s\, (u-s)^{2H}\rho (u-s).
\end{eqnarray*}
Let us first study $\beta_{\varepsilon }$. We can write
\[
\beta_{\varepsilon }= \e^{2H}\int_{0}^{T/\varepsilon }x^{2H}\rho (x)(T-\varepsilon x)
\,\mathrm{d}x.
\]
The integral $\int_{0}^{\infty }x^{2H}\rho (x)\,\mathrm{d}x$ is convergent for $H<1/4$,
while $\int_{0}^{T/\varepsilon }x^{2H}\rho (x)\,\mathrm{d}x$
diverges as $-\frac18\log(1/\e)$ for $H=1/4$ and
as $H(2H-1)T^{4H-1}\e^{1-4H}$ for $1/4<H<1/2$.
The integral $\int_{0}^{T/\varepsilon }x^{2H+1}\rho (x)\,\mathrm{d}x$ diverges as $H(2H-1)T^{4H}
\varepsilon^{-4H}$.
Therefore,
\[
\lim_{\varepsilon \rightarrow 0}\beta_{\varepsilon }=0.
\]
Second, let us write $\alpha_\e$ as
\begin{eqnarray*}
\alpha_{\varepsilon } &=& \varepsilon^{1+2H} \int_{0}^{T/\varepsilon
}\mathrm{d}u\, u^{2H}\int_{0}^{T/\e}\mathrm{d}s\, \rho (u-s)\\
&=& \varepsilon^{1+2H} \biggl( \int_{0}^{T/\varepsilon
}\mathrm{d}u\, u^{2H}\int_{0}^{u}\mathrm{d}s\, \rho (u-s) +  \int_{0}^{T/\varepsilon
}\mathrm{d}u\, u^{2H}\int_{u}^{T/\e}\mathrm{d}s\, \rho (u-s) \biggr)\\
&=& \frac{1}{2H+1} \int_{0}^{T/\varepsilon
}\rho(x) \bigl(T^{2H+1}-(\e x)^{2H+1}+(T-\e x)^{2H+1} \bigr)\,\mathrm{d}x.
\end{eqnarray*}
Hence, by the dominated convergence theorem, we have
\[
\lim_{\varepsilon \rightarrow 0}\alpha_{\varepsilon }=\frac{2T^{2H+1}}{2H+1}\int_0^\infty \rho(x)\,\mathrm{d}x
\]
so that
\[
\lim_{\varepsilon \rightarrow 0}\e^{1-2H}E[\widetilde{G}_{\varepsilon }^2]
=\lim_{\varepsilon \rightarrow 0}E[\widetilde{F}_{\varepsilon }^2]=
\frac{2T^{2H+1}}{2H+1}\int_0^\infty \rho(x)\,\mathrm{d}x.
\]

\textit{Step \textup{2}. Showing the convergence in law in} (\ref{star1}).

By the previous step, the distributions of the family
\[
\bigl( \bigl(B^{(1)}_t,B^{(2)}_t\bigr) _{t\in [0,T]},  \e^{\sfrac 12-H} \widetilde{G}_{\e}   \bigr)_{\e>0}
\]
are tight in $C([0,T]^2) \times \mathbb{R}$ and it suffices to show the convergence of the finite-dimensional distributions.
We need to show that  for any $\lambda \in \mathbb{R}$,
any $0<t_1\leq\cdots\leq t_k$, any $\theta_1,\ldots,\theta_k\in\R$ and
any $\mu_1,\ldots,\mu_k\in\R$,
we have
%e6.2 ###
\begin{eqnarray}\label{stardesstars}
&&\lim_{\varepsilon \downarrow 0}  E \bigl[  \mathrm{e}^{\mathrm{i}\sum_{j=1}^k\theta_jB_{t_j}^{(1)}}
\mathrm{e}^{\mathrm{i}\sum_{j=1}^k \mu_j B^{(2)}_{t_j}} \mathrm{e}^{\mathrm{i}\lambda \e^{ \sfrac 12 - H} \widetilde{G}_{\varepsilon }}  \bigr] \nonumber
\\[-8pt]\\[-8pt]
&&\quad=
E \bigl[
\mathrm{e}^{-(\sfrac12)\operatorname{Var} (\sum_{j=1}^k \mu_j B^{(2)}_{t_j} )}
\bigr]   E \bigl[  \mathrm{e}^{\mathrm{i}\sum_{j=1}^k\theta_jB_{t_j}^{(1)}}
\mathrm{e}^{- \sfrac {\lambda^2 S^2} 2}   \bigr],\nonumber
\end{eqnarray}
where  $S=\sqrt{2\int_0^\infty \rho(x)\,\mathrm{d}x\int_0^T( B^{(1)}_u)^2 \,\mathrm{d}u}$.
We can write
\begin{eqnarray*}
&&E \bigl[
\mathrm{e}^{\mathrm{i}\sum_{j=1}^k\theta_jB_{t_j}^{(1)}}
\mathrm{e}^{\mathrm{i}\sum_{j=1}^k \mu_j B^{(2)}_{t_j}} \mathrm{e}^{\mathrm{i}\lambda \e^{ \sfrac 12 - H} \widetilde{G}_{\varepsilon }}  \bigr]\\
&&\quad = E \bigl[  \mathrm{e}^{\mathrm{i}\sum_{j=1}^k\theta_jB_{t_j}^{(1)}} E\bigl[ \mathrm{e}^{\mathrm{i}\sum_{j=1}^k \mu_j B^{(2)}_{t_j}}
\mathrm{e}^{\mathrm{i}\lambda \e^{ \sfrac 12 - H} \widetilde{G}_{\varepsilon }} \big| B^{(1)}\bigr]   \bigr]\\
&&\quad = E \bigl[  \mathrm{e}^{\mathrm{i}\sum_{j=1}^k\theta_jB_{t_j}^{(1)}}
\mathrm{e}^{-\lambda \e^{\sfrac12-H}\sum_{j=1}^k\mu_j\int_0^T B^{(1)}_u E (B_{t_j}^{(2)}\times
\sfrac{B^{((2)}_{u+\e}-B^{(2)}_u)}{\e} )\,\mathrm{d}u}
\\
&&\qquad\hspace*{9pt}{}\times \mathrm{e}^{-(\sfrac{\lambda^2}2) \e^{1-2H}\int_{[0,T]^2} B^{(1)}_uB^{(1)}_v \rho_\e(u-v)\,\mathrm{d}u\,\mathrm{d}v}
\mathrm{e}^{-(\sfrac12)\operatorname{Var} (\sum_{j=1}^k \mu_j B^{(2)}_{t_j} )}\bigr]
\end{eqnarray*}
with $\rho_\e(x)=\frac12 (|x+\e|^{2H}+|x-\e|^{2H}-2|x|^{2H} )$.
Observe that
\[
\int_{[0,T]^2} B^{(1)}_uB^{(1)}_v\rho_\e(u-v)\,\mathrm{d}u\,\mathrm{d}v\geq 0
\]
since $\rho_\e(u-v)=E [(B^{(2)}_{u+\e}-B^{(2)}_u)(B^{(2)}_{v+\e}-B^{(2)}_v) ]$ is
a covariance function.
Moreover, for any fixed $t\geq 0$, we have
\begin{eqnarray*}
&&\int_0^T B^{(1)}_u E \biggl(B_{t}^{(2)}\times
\frac{B^{(2)}_{u+\e}-B^{(2)}_u}{\e} \biggr)
\,\mathrm{d}u\\
&&\quad=
\frac1{2}\int_0^T B^{(1)}_u
\biggl(\frac{(u+\e)^{2H}-u^{2H}}{\e}+\frac{|t-u|^{2H}-|t-u-\e|^{2H}}{\e}\biggr)\,\mathrm{d}u\\
&&\quad\stackrel{\mathrm{a.s.}}{\mathop{\longrightarrow}_{\e\to 0}}
H\int_0^T B^{(1)}_u   (u^{2H-1}-|t-u|^{2H-1} )\,\mathrm{d}u.
\end{eqnarray*}
Since $H<1/2$, this implies that
\[
\mathrm{e}^{-\lambda \e^{\sfrac12-H}\sum_{j=1}^k\mu_j\int_0^T B^{(1)}_u E
(B_{t_j}^{(2)}\times \afrac{B^{(2)}_{u+\e}-B^{(2)}_u}{\e} )
\,\mathrm{d}u}  \stackrel{\mathrm{a.s.}}{\mathop{\longrightarrow}_{\e\to 0}}   1.
\]
Hence, to get (\ref{stardesstars}), it suffices to show that
%e6.3 ###
\begin{equation}\label{stardes}
E
\bigl[
\mathrm{e}^{\mathrm{i}\sum_{j=1}^k\theta_jB_{t_j}^{(1)}}
\mathrm{e}^{-(\sfrac{\lambda^2}2) \e^{1-2H}\int_{[0,T]^2} B^{(1)}_uB^{(1)}_v \rho_\e(u-v)\,\mathrm{d}u\,\mathrm{d}v}
\bigr]
\mathop{\longrightarrow}_{\e\to 0}E \bigl[
\mathrm{e}^{\mathrm{i}\sum_{j=1}^k\theta_jB_{t_j}^{(1)}} \mathrm{e}^{-(\sfrac{\lambda^2}2)S^2} \bigr].
\end{equation}
We have
\begin{eqnarray*}
C_\e&:=&E
\Biggl[\exp
\Biggl(\mathrm{i}\sum_{j=1}^k\theta_jB_{t_j}^{(1)}
-\frac{\lambda^2}2 \e^{1-2H}\int_{[0,T]^2} B^{(1)}_uB^{(1)}_v \rho_\e(u-v)\,\mathrm{d}u\,\mathrm{d}v \Biggr)
\Biggr]\\
&=&E
\Biggl[\exp
\Biggl(\mathrm{i}\sum_{j=1}^k\theta_jB_{t_j}^{(1)}
-\lambda^2 \e^{1-2H}\int_0^T B^{(1)}_u \biggl(\int_0^u B^{(1)}_{u-x} \rho_\e(x)\,\mathrm{d}x \biggr)\,\mathrm{d}u \Biggr)
\Biggr]\\
&=&E
\Biggl[\exp
\Biggl(\mathrm{i}\sum_{j=1}^k\theta_jB_{t_j}^{(1)}
-\lambda^2 \e^{1-2H}\int_0^T \rho_\e(x) \biggl(\int_x^T B^{(1)}_uB^{(1)}_{u-x} \,\mathrm{d}u \biggr)\,\mathrm{d}x \Biggr)
\Biggr]\\
&=&E
\Biggl[\exp
\Biggl(\mathrm{i}\sum_{j=1}^k\theta_jB_{t_j}^{(1)}
-\lambda^2 \int_0^{T/\e} \rho(x) \biggl(\int_{\e x}^T B^{(1)}_uB^{(1)}_{u-\e x} \,\mathrm{d}u \biggr)\,\mathrm{d}x \Biggr)
\Biggr],
\end{eqnarray*}
the last inequality following from the relation $\rho_\e(x)=\e^{2H}\rho(x/\e)$.
By the dominated convergence theorem, we obtain
\begin{eqnarray*}
C_\e&\displaystyle\mathop{\longrightarrow}_{\e\to 0} &
E
\Biggl[\exp
\Biggl(\mathrm{i}\sum_{j=1}^k\theta_jB_{t_j}^{(1)}
-\lambda^2 \int_0^{\infty} \rho(x)\,\mathrm{d}x\times \int_{0}^T \bigl(B^{(1)}_u\bigr)^2 \,\mathrm{d}u \Biggr)
\Biggr]\\
&=&E
\Biggl[\exp
\Biggl(\mathrm{i}\sum_{j=1}^k\theta_jB_{t_j}^{(1)}
-\frac{\lambda^2}2 S^2 \Biggr)
\Biggr],
\end{eqnarray*}
that is, (\ref{stardes}). The proof of the theorem is thus completed.
\end{pf}

Recall the definition (\ref{rho}) of $\rho$ and the definition of $\breve{G}_\e$:
\[
\breve{G}_\e=
\int_0^T  \biggl(\int_0^u \frac{B^{(1)}_{v+\e}-B^{(1)}_v}{\e} \,\mathrm{d}v \biggr)\frac{B^{(2)}_{u+\e}-B^{(2)}_u}{\e} \,\mathrm{d}u.
\]
%t13
\begin{theo}\label{cvgausslevyareabis}
Convergences in law (\ref{star10}) and (\ref{star10crit}) hold.
\end{theo}

\begin{pf}
We only show the first convergence, the proof of the second one being very similar.
By using the scaling properties
of the fBm, first observe  that $\varepsilon ^{\sfrac{1}{2}-2H}\breve{G}_\e$ has the same law
as
\[
\breve{F}_{\varepsilon }=\sqrt{\e} \int_{0}^{T/\varepsilon}
\biggl(\int_0^u
\bigl(B_{v+1}^{(1)}-B^{(1)}_v \bigr)\,\mathrm{d}v\biggr)
\bigl( B_{u+1}^{(2)}-B_{u}^{(2)} \bigr) \,\mathrm{d}u.
\]

We now fix $H<1/4$ and the proof is divided into several steps.

\textit{Step \textup{1}. Computing the variance of  $\breve{F}_{\varepsilon }$}.
We can write
\begin{eqnarray*}
E(\breve{F}_{\varepsilon }^{2}) &=&\varepsilon \int_{[0,T/\varepsilon]^2}
\mathrm{d}u\,\mathrm{d}u' \,\rho(u-u' )\int_{0}^{u }\mathrm{d}v\int_0^{u' }\mathrm{d}v'\,\rho(v-v')
\end{eqnarray*}
with $\rho(x)=\frac12 (|x+1|^{2H}+|x-1|^{2H}-2|x|^{2H} )$.
We have
\[
\int_0^u\mathrm{d}v\int_0^{u' }\mathrm{d}v'\,\rho(v-v')=\frac{\Psi(u-u' )-\Psi(u)-\Psi(u'
)+2}{2(2H+1)(2H+2)},
\]
where
%e6.4 ###
\begin{equation}\label{PPSSII}
\Psi(x)=2|x|^{2H+2}-|x+1|^{2H+2}-|x-1|^{2H+2}.
\end{equation}
Consider first the contribution of the term $\Psi(u-u' )$. We have
\[
\lim_{\e\to 0}\e\int_{[0,T/\e]^2}\rho(u-u' )\Psi(u-u' )\,\mathrm{d}u\,\mathrm{d}u' =T\int_\R\rho(x)\Psi(x)\,\mathrm{d}x.
\]
Note that $\rho(x)\sim H(2H-1)|x|^{2H-2}$ and
$\Psi(x)\sim -(2H+2)(2H+1)|x|^{2H}$ as $|x|\to\infty$
so that $\int_\R|\rho(x)\Psi(x)|\,\mathrm{d}x<\infty$ because $H<1/4$. On the other hand,
we have
\[
\e\int_{[0,T/\e]^2} \rho(u-u' )\Psi(u)\,\mathrm{d}u\,\mathrm{d}u'
=\e\int_0^{T/\e}\mathrm{d}u\,\Psi(u)\int_{u-T/\e}^u\mathrm{d}x\,\rho(x)
\]
and this converges to zero as $\e\to 0$. Indeed,
since $\rho(x)\sim H(2H-1)x^{2H-2}$ as $x\to\infty$,
we have $\int_u^\infty \rho(x)\,\mathrm{d}x \sim Hu^{2H-1}$ as $u\to\infty$;
hence,
since $\int_\R\rho(x)\,\mathrm{d}x=0$, $H<1/4$ and
$\Psi(u)\sim -(2H+2)(2H+1)u^{2H}$ as $u\to\infty$, we have
\[
\lim_{u\to\infty}\Psi(u)\int_{-\infty}^u\rho(x)\,\mathrm{d}x
=-\lim_{u\to\infty}\Psi(u)\int_{u}^\infty\rho(x)\,\mathrm{d}x
=0.
\]
Also, we have
\[
\lim_{\e\to 0}\e\int_{[0,T/\e]^2}\rho(u-u' )\,\mathrm{d}u\,\mathrm{d}u' =\int_\R\rho(x)\,\mathrm{d}x=0.
\]
Therefore, $\lim_{\e\to 0}E(\breve{F}_\e^2)=\breve{\sigma}_H^2$.

\textit{Step \textup{2}. Showing the convergence in law} (\ref{star10}).
We first remark that by step 1, the laws of the family
$ ((B^{(1)}_t,B^{(2)}_t) _{t\in [0,T]},\e^{\sfrac12-2H}\breve{G}_{\varepsilon } ) _{\e>0}$ are tight.
Therefore, we only have to prove the convergence of the finite-dimensional laws.
Moreover, by the main result of Peccati and Tudor \cite{PT}, it suffices to prove  that
%e6.5 ###
\begin{equation}\label{felow}
\e^{\sfrac12-2H}\breve{G}_\e\stackrel{\mathrm{Law}}{=}\breve{F}_{\varepsilon }\stackrel{\mathrm{Law}}{\longrightarrow}
\mathscr{N}(0,T\breve{\sigma}^2_H) \qquad\mbox{as $\e\to 0$}.
\end{equation}
We have
\begin{eqnarray*}
E(\mathrm{e}^{\mathrm{i}\lambda \breve{F}_\e})
&=&E
\biggl(\exp \biggl\{-\frac{\lambda^2\e}2\int_{\lbrack 0,T/\varepsilon]^{2}}
\bigl(B^{(2)}_{u+1}-B^{(2)}_{u}\bigr)\bigl(B^{(2)}_{u' +1}-B^{(2)}_{u' }\bigr)
\\
&&\hspace*{93pt}{}\times
\biggl( \int_{0}^{u}\int_{0}^{u'}\rho
(v-v^{\prime })\,\mathrm{d}v\,\mathrm{d}v^{\prime } \biggr) \,\mathrm{d}u\,\mathrm{d}u'  \biggr\}
\biggr).
\end{eqnarray*}
Since $\rho(v-v')=E [(B^{(1)}_{v+1}-B^{(1)}_v)(B^{(1)}_{v'+1}-B^{(1)}_v) ]$ is
a covariance function, observe that the quantity inside the exponential in the right-hand side of the previous identity is
negative.
Hence, since $x\mapsto\exp (-\frac{\lambda^2}2x_+ )$
is continuous and bounded by 1 on $\R$, (\ref{felow}) will be a consequence of the  convergence
%e6.6 ###
\begin{equation}\label{toshow2}
A_\e \stackrel{\rm law}{\longrightarrow}T\breve{\sigma}_H^2 \qquad\mbox{as $\e\to 0$}
\end{equation}
with
\[
A_{\varepsilon}:=\varepsilon \int_{\lbrack 0,T/\varepsilon
]^{2}}(B_{u+1}-B_{u})(B_{u' +1}-B_{u' }) \biggl( \int_{0}^{u }\int_{0}^{u'}\rho
(v-v^{\prime })\,\mathrm{d}v\,\mathrm{d}v^{\prime } \biggr) \,\mathrm{d}u\,\mathrm{d}u',
\]
$B$ denoting a fractional Brownian motion with Hurst index $H$.
The proof of (\ref{toshow2}) will be achieved by showing that the expectation (resp., the variance)
of $A_\e$ tends to $T\breve{\sigma}_H^2$ (resp., zero).
By step~1, observe that
\[
E(A_\e)=E(\breve{F}_\e^2)\to T\breve{\sigma}^2_H
\]
as $\e\to 0$. We now want to show that the variance of $A_{\varepsilon }$ converges to zero.
Performing the changes of variables $s=u\varepsilon$ and $t=u' \varepsilon$ yields
\[
A_{\varepsilon }=\varepsilon ^{-1}\int_{[0,T]^{2}}(B_{s/\varepsilon
+1}-B_{s/\varepsilon })(B_{t/\varepsilon +1}-B_{t/\varepsilon})
\biggl(\int_{0}^{s/\varepsilon }\int_{0}^{t/\varepsilon }\rho (v-v^{\prime })\,\mathrm{d}v\,\mathrm{d}v^{\prime } \biggr) \,\mathrm{d}s\,\mathrm{d}t,
\]
which has the same distribution as
\begin{eqnarray*}
C_{\varepsilon } &=&\varepsilon ^{-1-2H}\int_{[0,T]^{2}}(B_{s+\varepsilon
}-B_{s})(B_{t+\varepsilon }-B_{t}) \biggl( \int_{0}^{s/\varepsilon
}\int_{0}^{t/\varepsilon }\rho (u-u' )\,\mathrm{d}u\,\mathrm{d}u'  \biggr) \,\mathrm{d}s\,\mathrm{d}t \\
&=&\varepsilon ^{-1-2H}\int_{[0,T]^{2}}(B_{s+\varepsilon
}-B_{s})(B_{t+\varepsilon }-B_{t})\Lambda _{\varepsilon }(s,t)\,\mathrm{d}s\,\mathrm{d}t,
\end{eqnarray*}
where $\Lambda _{\varepsilon }(s,t)=\int_{0}^{s/\varepsilon
}\int_{0}^{t/\varepsilon }\rho (u-u' )\,\mathrm{d}u\,\mathrm{d}u'.$  This can be written as
\[
C_{\varepsilon }=\varepsilon ^{-1-2H}\int_{\mathbb{R}^{2}}B_{s}B_{t}\Sigma_{\varepsilon }(s,t)\,\mathrm{d}s\,\mathrm{d}t,
\]
where
%e6.7 ###
\begin{eqnarray}\label{e1}
\hspace*{-15pt}\Sigma _{\varepsilon }(s,t) &=&\mathbf{1}_{[\varepsilon,T+\varepsilon ]}(s)
\mathbf{1}_{[\varepsilon,T+\varepsilon ]}(t)\Lambda _{\varepsilon}(s-\varepsilon,t-\varepsilon )
-\mathbf{1}_{[0,T]}(s)\mathbf{1}_{[\varepsilon,T+\varepsilon ]}(t)\Lambda _{\varepsilon }(s,t-\varepsilon )
\nonumber \\[-8pt]\\[-8pt]
&&{}-\mathbf{1}_{[\varepsilon,T+\varepsilon ]}(s)\mathbf{1}_{[0,T]}(t)\Lambda_{\varepsilon }
(s-\varepsilon,t)+\mathbf{1}_{[0,T]}(s)\mathbf{1}_{[0,T]}(t)\Lambda _{\varepsilon }(s,t).\nonumber
\end{eqnarray}
Moreover,
\[
C_{\varepsilon }-E(C_{\varepsilon })=\varepsilon ^{-1-2H}I_2 \biggl(\int_{\mathbb{R}^{2}}
\mathbf{1}_{[0,s]}\otimes \mathbf{1}_{[0,t]} \Sigma_{\varepsilon }(s,t)\,\mathrm{d}s\,\mathrm{d}t \biggr),
\]
where $I_{2}$ is the double stochastic integral with respect to $B$.
Therefore,
\begin{eqnarray*}
\operatorname{Var}(C_{\varepsilon }) &=&2\varepsilon ^{-2-4H}
\bigg\| \int_{\mathbb{R}^{2}} \mathbf{1}_{[0,s]}\otimes \mathbf{1}_{[0,t]}\Sigma _{\varepsilon }(s,t)\,\mathrm{d}s\,\mathrm{d}t\bigg\|^2_{\HH^{\otimes 2}} \\
&=&2\varepsilon ^{-2-4H}\int_{\mathbb{R}^{4}}R_{H}(s,s^{\prime
})R_{H}(t,t^{\prime })\Sigma _{\varepsilon }(s,t)\Sigma _{\varepsilon
}(s^{\prime },t^{\prime })\,\mathrm{d}s\,\mathrm{d}t\,\mathrm{d}s^{\prime }\,\mathrm{d}t^{\prime }.
\end{eqnarray*}
Taking into account that the partial derivatives $\frac{\partial R_{H}}{\partial s}$ and
$\frac{\partial R_{H}}{\partial t}$
are integrable, we can write
\begin{eqnarray*}
&& \operatorname{Var}(C_{\varepsilon })=2\varepsilon ^{-2-4H}\int_{\mathbb{R}^{4}}
\biggl( \int_{0}^{s}\frac{\partial R_{H}}{\partial\sigma }(\sigma,s^{\prime})\,\mathrm{d}\sigma  \biggr)
\biggl( \int_{0}^{t^{\prime }}\frac{\partial R_{H}}{\partial \tau }(t,\tau )\,\mathrm{d}\tau  \biggr)   \\
&&\hspace*{96pt} {}    \times  \Sigma _{\varepsilon }(s,t)\Sigma_{\varepsilon }(s^{\prime },t^{\prime })
\,\mathrm{d}s\,\mathrm{d}t\,\mathrm{d}s^{\prime }\,\mathrm{d}t^{\prime}.
\end{eqnarray*}
Hence, by integrating by parts, we get
\begin{eqnarray*}
&& \operatorname{Var}(C_{\varepsilon })=2\varepsilon ^{-2-4H}\int_{\mathbb{R}^{4}}\frac{\partial R_{H}}{\partial s}(s,s^{\prime })
\frac{\partial R_{H}}{\partial t^{\prime }}(t,t^{\prime })    \\
&&\hspace*{97pt}{}   \times  \biggl( \int_{0}^{s}
\times \Sigma _{\varepsilon
}(\sigma,t)\,\mathrm{d}\sigma  \biggr)  \biggl( \int_{0}^{t^{\prime }}\Sigma_{\varepsilon }(s^{\prime },\tau )\,\mathrm{d}\tau  \biggr)
\,\mathrm{d}s\,\mathrm{d}t\,\mathrm{d}s^{\prime }\,\mathrm{d}t^{\prime}.
\end{eqnarray*}
From (\ref{e1}), we obtain
\[
\int_{0}^{s}\Sigma _{\varepsilon }(\sigma,t)\,\mathrm{d}\sigma =\mathbf{1}_{[0,T]}(s) \bigl( \mathbf{1}_{[0,\varepsilon ]}(t)-
\mathbf{1}_{[T,T+\varepsilon ]}(t) \bigr) \int_{s-\varepsilon }^{s}\Lambda_{\varepsilon }(\sigma,t-\varepsilon )\,\mathrm{d}\sigma.
\]
In the same way,
\[
\int_{0}^{t^{\prime }}\Sigma _{\varepsilon }(s^{\prime },\tau )\,\mathrm{d}\tau =
\mathbf{1}_{[0,T]}(t^{\prime }) \bigl( \mathbf{1}_{[0,\varepsilon
]}(s^{\prime })-\mathbf{1}_{[T,T+\varepsilon ]}(s^{\prime }) \bigr)
\int_{t^{\prime }-\varepsilon }^{t^{\prime }}\Lambda _{\varepsilon
}(s^{\prime }-\varepsilon,\tau )\,\mathrm{d}\tau.
\]
As a consequence,
\begin{eqnarray*}
\operatorname{Var}(C_{\varepsilon }) &=&2\varepsilon ^{-2-4H}\int_{\mathbb{R}^{4}}
\frac{\partial R_{H}}{\partial s}(s,s^{\prime })\frac{\partial R_{H}}
{\partial t^{\prime }}(t,t^{\prime }) \biggl( \int_{s-\varepsilon }^{s}\Lambda_{\varepsilon }(\sigma,t-\varepsilon )\,\mathrm{d}\sigma  \biggr)  \\
&&\hspace*{53pt}{}\times  \biggl( \int_{t^{\prime }-\varepsilon }^{t^{\prime }}
\Lambda_{\varepsilon }(s^{\prime }-\varepsilon,\tau )\,\mathrm{d}\tau  \biggr) \mathbf{1}_{[0,T]}(s)
\bigl( \mathbf{1}_{[0,\varepsilon ]}(t)-\mathbf{1}_{[T,T+\varepsilon ]}(t) \bigr)  \\
&&\hspace*{53pt}{}\times \mathbf{1}_{[0,T]}(t^{\prime }) \bigl( \mathbf{1}_{[0,\varepsilon]}(s^{\prime })-\mathbf{1}_{[T,T+\varepsilon ]}(s^{\prime }) \bigr)
\,\mathrm{d}s\,\mathrm{d}t\,\mathrm{d}s^{\prime }\,\mathrm{d}t^{\prime }
=\sum_{i=1}^{4}H_{\varepsilon }^{i},
\end{eqnarray*}
where
\begin{eqnarray*}
H_{\varepsilon }^{1} &=&\int_{0}^{T}\int_{0}^{\varepsilon
}\int_{0}^{\varepsilon }\int_{0}^{T}G_{\varepsilon }(s,t,s^{\prime
},t^{\prime })\,\mathrm{d}s\,\mathrm{d}t\,\mathrm{d}s^{\prime }\,\mathrm{d}t^{\prime}, \\
H_{\varepsilon }^{2} &=&-\int_{0}^{T}\int_{T}^{T+\varepsilon
}\int_{0}^{\varepsilon }\int_{0}^{T}G_{\varepsilon }(s,t,s^{\prime
},t^{\prime })\,\mathrm{d}s\,\mathrm{d}t\,\mathrm{d}s^{\prime }\,\mathrm{d}t^{\prime}, \\
H_{\varepsilon }^{3} &=&-\int_{0}^{T}\int_{0}^{\varepsilon
}\int_{T}^{T+\varepsilon }\int_{0}^{T}G_{\varepsilon }(s,t,s^{\prime
},t^{\prime })\,\mathrm{d}s\,\mathrm{d}t\,\mathrm{d}s^{\prime }\,\mathrm{d}t^{\prime}, \\
H_{\varepsilon }^{4} &=&\int_{0}^{T}\int_{0}^{T+\varepsilon
}\int_{0}^{T+\varepsilon }\int_{0}^{T}G_{\varepsilon }(s,t,s^{\prime
},t^{\prime })\,\mathrm{d}s\,\mathrm{d}t\,\mathrm{d}s^{\prime }\,\mathrm{d}t^{\prime}
\end{eqnarray*}
and
\begin{eqnarray*}
 G_{\varepsilon }(s,t,s^{\prime },t^{\prime })&=&2\varepsilon ^{-2-4H}
\frac{\partial R_{H}}{\partial s}(s,s^{\prime })\frac{\partial R_{H}}{\partial
t^{\prime }}(t,t^{\prime })    \\
&&{} \times  \biggl( \int_{s-\varepsilon }^{s}\Lambda_{\varepsilon }(\sigma,t-\varepsilon )\,\mathrm{d}\sigma  \biggr)
\biggl(\int_{t^{\prime }-\varepsilon }^{t^{\prime }}\Lambda _{\varepsilon
}(s^{\prime }-\varepsilon,\tau )\,\mathrm{d}\tau  \biggr).
\end{eqnarray*}
We only consider the term $H_{\varepsilon }^{1}$ because the others can
be handled in the same way.
We have, with $\Psi$ given by (\ref{PPSSII}),
\[
\Lambda _{\varepsilon }(s,t)=\int_{0}^{s/\varepsilon
}\int_{0}^{t/\varepsilon }\rho (u-u')\,\mathrm{d}u\,\mathrm{d}u'=\frac{\Psi (\sfrac{(s-t)}{\varepsilon })
-\Psi (\sfrac{s}{\varepsilon })-\Psi (\sfrac{t}{\varepsilon })+2}{2(2H+1)(2H+2)}.
\]
Note that
\begin{eqnarray*}
\bigg| \Psi \biggl(\frac{s-t}{\varepsilon }\biggr) \bigg|  &\le &\varepsilon
^{-2H-2} | 2|s-t|^{2H+2}-|s-t+\varepsilon |^{2H+2}-|s-t-\varepsilon
|^{2H+2} |  \\
&\leq &C\varepsilon ^{-2H}
\end{eqnarray*}
for any $s,t\in \lbrack 0,T]$. Therefore, $ | \Lambda _{\varepsilon
}(s,t) | \leq C\varepsilon ^{-2H}$ and we obtain   the
estimate
\[
| G_{\varepsilon }(s,t,s^{\prime },t^{\prime }) | \leq
C\varepsilon ^{-8H} ( s^{2H-1}+|s-s^{\prime }|^{2H-1} )
(t^{\prime 2H-1}+|t-t^{\prime }|^{2H-1} ).
\]
As a consequence,
\begin{eqnarray*}
| H_{\varepsilon }^{1} |  &\leq &\int_{0}^{T}\int_{0}^{\varepsilon
}\int_{0}^{\varepsilon }\int_{0}^{T} | G_{\varepsilon }(s,t,s^{\prime
},t^{\prime }) | \,\mathrm{d}s\,\mathrm{d}t\,\mathrm{d}s^{\prime }\,\mathrm{d}t^{\prime } \\
&\leq &C\varepsilon ^{-8H}\int_{0}^{T}\int_{0}^{\varepsilon
}\int_{0}^{\varepsilon }\int_{0}^{T} ( s^{2H-1}+|s-s^{\prime}|^{2H-1} )   \\
&&{}\hspace*{96pt}  \times  ( t^{\prime 2H-1}+|t-t^{\prime }|^{2H-1} )
\,\mathrm{d}s\,\mathrm{d}t\,\mathrm{d}s^{\prime }\,\mathrm{d}t^{\prime } \\
&\leq &C\varepsilon ^{2-8H},
\end{eqnarray*}
which converges to zero because $H<\frac{1}{4}$.
\end{pf}

Recall the definition (\ref{hat}) of $\widehat{G}_\e$,
\[
\widehat{G}_\e=
\int_0^T \frac{B^{(1)}_{u+\e}-B^{(1)}_u}{\e}\times\frac{B^{(2)}_{u+\e}-B^{(2)}_u}{\e} \,\mathrm{d}u.
\]
We have the following result.
%t14
\begin{theo}\label{cvgausscovariation}
Convergences (\ref{star5}) and (\ref{star6}) hold.
\end{theo}

\begin{pf}
We use the same trick as in \cite{nourdinjfa}, Remark 1.3, point 4.
Let $\beta$ and $\widetilde{\beta}$ be two independent one-dimensional fractional Brownian motions with index
$H$. Set $B^{(1)}=(\beta+\widetilde{\beta})/\sqrt{2}$ and $B^{(2)}=(\beta-\widetilde{\beta})/\sqrt{2}$.
It is easily checked that $B^{(1)}$ and $B^{(2)}$ are also two independent fractional Brownian
motions with index $H$. Moreover, we have
%e6.8 ###
\begin{eqnarray}\label{dollar}
\hspace*{-25pt}\e^{\sfrac32-2H}\widehat{G}_\e
&=&\frac12 \e^{\sfrac32-2H}\int_0^T  \biggl(\frac{\beta_{u+\e}-\beta_u}{\e} \biggr)^2\,\mathrm{d}u
-\frac12 \e^{\sfrac32-2H}\int_0^T  \biggl(\frac{\widetilde{\beta}_{u+\e}-\widetilde{\beta}_u}{\e} \biggr)^2\,\mathrm{d}u\nonumber\\
&=&\frac1{2\sqrt{\e}}\int_0^T  \biggl(\frac{\beta_{u+\e}-\beta_u}{\e^H} \biggr)^2\,\mathrm{d}u
-\frac1{2\sqrt{\e}}\int_0^T  \biggl(\frac{\widetilde{\beta}_{u+\e}-\widetilde{\beta}_u}{\e^H} \biggr)^2\,\mathrm{d}u\\
&=&\frac1{2\sqrt{\e}}\int_0^T h_2 \biggl(\frac{\beta_{u+\e}-\beta_u}{\e^H} \biggr)\,\mathrm{d}u
-\frac1{2\sqrt{\e}}\int_0^T h_2 \biggl(\frac{\widetilde{\beta}_{u+\e}-\widetilde{\beta}_u}{\e^H} \biggr)\,\mathrm{d}u.\nonumber
\end{eqnarray}
The proofs of the desired convergences in law
are now direct consequences of
the convergence (\ref{cv<})    with $k=2$, taking into account that $\beta$
and $\widetilde{\beta}$ are independent.
\end{pf}

%r15
\begin{rem}
As a by-product of the decomposition (\ref{dollar}), and taking into account (\ref{cv>}) for $k=2$,
we get that $\int_0^T \dot{B}^{(1)}_u\diamond \dot{B}^{(2)}_u\, \mathrm{d}u$ and $ (Z_T^{(2)}-\widetilde{Z}_T^{(2)} )/2$
have the same law when $H>3/4$, where $\widetilde{Z}_T^{(2)}$ stands for an independent copy of the
Hermite random variable
$Z_T^{(2)}$.
\end{rem}

\section*{Acknowedgements}

The research of I. Nourdin was supported in part by the ANR project `Exploration des Chemins Rugueux'.
The research of D. Nualart was supported by NSF Grant DMS-0904538.

\printhistory

\end{document}